\documentclass[a4paper,11pt]{article}

\scrollmode
\usepackage{amsmath,amssymb,amsfonts,graphicx,color,fancyhdr}
\usepackage[latin1]{inputenc}
\usepackage[colorlinks=true, pdfstartview=FitV, linkcolor=blue, citecolor=blue, urlcolor=blue]{hyperref}

\newtheorem{thm}{Theorem}[section]
\newtheorem{prop}[thm]{Proposition}
\newtheorem{lem}[thm]{Lemma}
\newtheorem{df}[thm]{Definition}

\newtheorem{rem}[thm]{Remark}
\newtheorem{cor}[thm]{Corollary}

\def\be#1 {\begin{equation} \label{#1}}
\newcommand{\ee}{\end{equation}}
\def\dem {\noindent {\bf Proof : }}

\def\sqw{\hbox{\rlap{\leavevmode\raise.3ex\hbox{$\sqcap$}}$%
\sqcup$}}
\def\findem{\ifmmode\sqw\else{\ifhmode\unskip\fi\nobreak\hfil
\penalty50\hskip1em\null\nobreak\hfil\sqw
\parfillskip=0pt\finalhyphendemerits=0\endgraf}\fi}

\newcommand{\mb}{\medskip\noindent}
\newcommand{\gb}{\bigskip\noindent}
\newcommand{\R}{\mathbb R}

\newcommand{\N}{\mathbb N}

\newcommand{\C}{\mathbb C}

\newcommand{\A}{{\mathcal A}}
\newcommand{\K}{\kappa}

\def\Xint#1{\mathchoice
   {\XXint\displaystyle\textstyle{#1}}%
   {\XXint\textstyle\scriptstyle{#1}}%
   {\XXint\scriptstyle\scriptscriptstyle{#1}}%
   {\XXint\scriptscriptstyle\scriptscriptstyle{#1}}%
   \!\int}
\def\XXint#1#2#3{{\setbox0=\hbox{$#1{#2#3}{\int}$}
     \vcenter{\hbox{$#2#3$}}\kern-.5\wd0}}

\def\aver#1{\Xint-_{#1}}

\newcommand{\M}{{\mathcal M}}

\title{A $T(1)$-Theorem in relation to a semigroup of operators and applications to new paraproducts}
\author{Fr\'ed\'eric Bernicot
\\ CNRS - Laboratoire Paul Painlev\'e \\ Universit\'e
Lille 1\\59655 Villeneuve d'Ascq Cedex, France\\frederic.bernicot@math.univ-lille1.fr }
\date{March 31, 2010}

\begin{document}

\maketitle

\begin{abstract} In this work, we are interested to develop new directions of the famous $T(1)$-theorem. More precisely, we develop a general framework where we look for replacing the John-Nirenberg space $BMO$ (in the classical result) by a new $BMO_{L}$, associated to a semigroup of operators $(e^{-tL})_{t>0}$. These new spaces $BMO_L$ (including $BMO$) have recently appeared in numerous works in order to extend the theory of Hardy and BMO space to more general situations. Then we give applications by describing boundedness for a new kind of paraproducts, built on the considered semigroup. In addition we obtain a version of the classical T(1) theorem for doubling Riemannian manifolds.
\end{abstract}

\mb {\bf Key-words:} T(1)-theorem ; semigroup of operators ; paraproducts.

\mb {\bf MSC:} 30E20 ; 42B20 ; 42B30.

\tableofcontents

\section{Introduction}

The $T(1)$ Theorem of G. David and J. L. Journ\'e provides a very powerful tool for analyzing the $L^2$-boundedness of a Calder\'on-Zygmund. It claims that for such a linear operator $T$ if it satisfies a weak boundedness property then $T$ is bounded on $L^2$ if and only if $T(1)$ and $T^*(1)$ belong to the John-Nirenberg space $BMO$ (introduced in \cite{JN}). This result proved in the Euclidean space framework in \cite{DJ}, is based on an appropriate reproducing Calder\'on formula and the notion of Carleson measure (closely related to the space $BMO$). It permits for example to obtain a new proof of the $L^2$-boundedness of the first Calder\'on commutator and to obtain a characterization of the $L^2$-bounded pseudo-differential operators (belonging to some exotic class like $S^0_{1,1}$, see \cite{cm}). \\
Then this result was extended by F. Nazarov, S. Treil and A. Volberg for a non-doubling measure in \cite{NTV} and by X. Tolsa in proper subset of the Euclidean space (not satisfying the doubling property, see \cite{T}). Then many works deal with the following problem: change the function $1$ by an accretive function $b$ getting the so-called ``T(b) Theorem'' (see the work of A. McIntosh and Y. Meyer \cite{McM} and of G. David, J.L Journ\'e ans S. Semmes \cite{DJS})
 or by a system of accretive functions, obtaining the so-called ``local T(b) Theorem'' (see the work of M. Christ \cite{Ch} and other works \cite{NTV2}, ...). We refer the reader to a survey of S. Hofmann \cite{Ho} about such questions and related applications to PDEs.
Numerous works deal with some adaptations, for example to quadratic T(1) type theorem (see \cite{AT}), to off-diagonal T(1) theorem (see \cite{H}), to Triebel-Lizorkin spaces (see \cite{Y} and \cite{V}) or to vector-valued opertors (see \cite{Hy, HW}). 

\gb

All these results concerns Calder\'on-Zygmund operators and involve the John-Nirenberg space BMO. This space naturally arises as the class of functions whose deviation from their means over cubes is bounded. This space is strictly including the $L^\infty$ space and is a good extension of the Lebesgue spaces scale $(L^p)_{1<p<\infty}$ for $p\to \infty$ from a point of view of Harmonic Analysis. For example, it plays an important role in boundedness of Calder\'on-Zygmund operators, real interpolation, Carleson measure, study of paraproducts,  ... \\
 Unfortunately, there are situations where the John-Nirenberg space $BMO$ is not the
 right substitute to $L^\infty$ and there have been recently numerous works
 whose the goal is to define an adapted BMO space according to the context (see \cite{DY, DY1, HM} ...).
 For example the classical space $BMO$ is not well adapted to operators such as the Riesz transform
  on Riemannian manifolds. That is why in \cite{HM}, S. Hofmann and S. Mayboroda develop theory of Hardy and BMO spaces associated to a second order divergence form elliptic operators, which also including the corresponding John-Nirenberg inequality. In the recent works \cite{DY} and \cite{DY1}, X. T. Duong and L. Yan studied some new BMO type spaces and proved an associated version of John-Nirenberg inequality on these spaces (with duality results). In \cite{BZ,BZ3}, J. Zhao and the author have developed an abstract framework for BMO spaces (and proved some results about John-Nirenberg inequalities). 
This framework permits to cover the classical space $BMO$ and those defined in \cite{DY} and \cite{HM}. 

\gb The aim of this article is to continue the study of T(1) theorems using these new spaces BMO, defined by a semigroup of operators.

\gb 

Let us just describe some motivation to such results. Associated to the $BMO$ space, which is a good extension of the Lebesgue spaces scale $(L^p)_{1<p<\infty}$ for $p\to \infty$, there is the ``Hardy spaces'' $H^1$ for $p\to 1$.
On the Euclidean space, R. Coifman and G. Weiss have introduced the first Hardy space via atomic decomposition in \cite{CW}. Then several characterizations was obtained in \cite{stein}, via maximal function or Riesz transform. 
The equivalence between all these definitions due to R. Coifman in \cite{FS} can be understood from the celebrated theorem of C. Fefferman which says 
$$ ({H}^1)^{*} = BMO.$$
As $BMO$ for $L^\infty$, the space $H^1$ is a good substitute of $L^1(X)$ for many reasons. For instance, Calder\'on-Zygmund operators map $H^1$ to $L^1$. In addition, $H^1$ (and its dual) interpolates with Lebesgue spaces $L^p$, $1<p<\infty$. 
However, there are situations where the space $H^1$ is not the right substitute to $L^1$ and there has been recently a number of works with goal to define an adapted Hardy space \cite{A1,A,AR,BZ,B,DY,DY1,D1,D,D2,DZ,HM}.
In \cite{BZ, B}, the authors have described a very abstract theory for Hardy spaces (built via atomic decomposition) and interpolation results with Lebesgue spaces. This part seems to be well understand. Specially concerning Hardy spaces associated to semigroup, we have several characterizations via atomic decomposition, area square function and maximal function. \\
In order to pursue this theory, it is now important to get criterions for an operator to be bounded in $L^2$. Then using this theory of Hardy space, we know how to obtain $L^p$-boundedness from the initial $L^2$-boundedness. Aiming that, we are motivated to obtain a general T(1) theorem (associated to this framework of semigroup).

\gb
So consider an operator $L$ of order $m$ acting on a doubling Riemannian manifold $(M,d,\mu)$, such that it admits an holomorphic calculus. In this case, we can consider the semigroup $(e^{-tL})_{t>0}$. We define the space $BMO_L$ as the set of functions $f$ such that
$$ \sup_{t>0} \ \sup_{\genfrac{}{}{0pt}{}{Q \textrm{ball}}{r_Q^m=t}} \ \frac{1}{\mu(Q)}\int_Q \left|f-e^{-tL}f\right| d\mu <\infty.$$
We refer the reader to Section \ref{sec:pre} for precise assumptions on the manifold and on the semigroup and preliminaries on this new BMO space.
Then our main theorem is the following one (see Theorem \ref{thm:principal})~:

\begin{thm} \label{thm:thm} Let $T$ be a linear operator, weakly continuous on $L^2(M)$ and admitting off-diagonal decays relatively to cancellation built with the semigroup (see Assumptions (\ref{assum:off}), (\ref{assum:off-dual}) and   (\ref{assum:weaklybounded})). Then if $T(1)\in BMO_L$ and $T^*(1)\in BMO_{L^*}$ then $T$ admits a bounded exstension in $L^2(M)$. 
\end{thm}
Moreover we will describe a reverse property: if $T$ admits a $L^2$-bounded extension then $T(1)$ and $T^*(1)$ belong to some BMO spaces (closely related to $BMO_L$ and $BMO_{L^*}$). 

\mb So we obtain as for the original theorem of G. David and J.L. Journ\'e, a criterion for the $L^2$-boundedness via these new spaces $BMO_L$. We emphasize that Assumptions (\ref{assum:off}), (\ref{assum:off-dual}) and   (\ref{assum:weaklybounded}) can be just considered as generalizations of usual Calder\'on-Zygmund properties since if we consider $L=-\Delta$ on $\R^d$ then every Calder\'on-Zygmund operators satisfy to (\ref{assum:off}), (\ref{assum:off-dual}) and (\ref{assum:weaklybounded}). So we recover the original T(1) theorem (since $BMO_{L}=BMO$ in this case).

\mb Then in Section \ref{sec:para}, we apply this new T(1) theorem to a new kind of paraproducts. In the usual framework, paraproducts were introduced by J. M. Bony \cite{bony} and then studied by Y. Meyer in \cite{meyer1, meyer2}. These operators are considered as the prototype of the so called ``Coifman-Meyer bilinear operators'', which are particular bilinear Calder\'on-Zygmund. A main result says that they are bounded from $L^p \times L^q$ into $L^{r'}$ as soon as 
$1<p,q\leq \infty$ and
$$ \frac{1}{p}+\frac{1}{q}=\frac{1}{r'}.$$
These operators appear in the study of the pointwise product between two functions and their study require to understand the frequency analysis of the product. Here, we introduce new bilinear operators, built with semigroups, which correspond to paraproducts with a frequency analysis adapted to the spectral properties of our semigroups $e^{-tL}$. In this case, it is not clear how $e^{-tL}$ acts on a product of two functions and that is why boundedness of such operators are not clear (see Remark \ref{rem:r}). We refer the reader to Section \ref{sec:para} for precise statement of the results. We prove the following one:

\begin{thm}[Theorem \ref{thm:paraproduit3}] Let $\psi$ and $\phi$ be defined as in Theorems \ref{thm:paraproduit} and \ref{thm:paraproduit2}.
Let $p,q\in (1,\infty)$ and $r'\in(1/2,\infty)$ be exponents satisfying
$$ \frac{1}{p}+\frac{1}{q} = \frac{1}{r'}$$
The paraproducts
$$ (f,h) \rightarrow \int_{0}^\infty \psi_t(L)\left[\phi_t(L)f \, \phi_t(L)h\right] \frac{dt}{t}$$
and
$$ (h,f) \rightarrow \int_{0}^\infty \phi_t(L)\left[\psi_t(L)f \, \phi_t(L)h\right] \frac{dt}{t} $$
are bounded from $L^p(M) \times L^q(M)$ to $L^{r'}(M)$.
\end{thm}
Moreover we obtain weighted estimates.

\mb We finish this work by describing in Subsection \ref{subsec:manifold}, a version of the classical T(1) theorem for Calder\'on-Zygmund operators in general doubling Riemannian manifold. All the previous cited works only deal with the Euclidean space and are based on some specific differential properties of $\R^n$. That is why, it seems important to us to prove such result, and we obtain it as an application of our main result.

\section{Preliminaries} \label{sec:pre}

 For a ball $Q$ in a metric space, $\lambda Q$  denotes the ball co-centered with $Q$ and with radius $\lambda$ times that of $Q$. Finally, $C$ will be a constant that may change from an inequality to another and we will use $u\lesssim
v$ to say that there exists a constant $C$  such that $u\leq Cv$ and $u\simeq v$ to say that $u\lesssim v$ and $v\lesssim u$.

\mb In all this paper, $M$ denotes a complete Riemannian manifold. We write $\mu$ for the Riemannian measure on $M$, $\nabla$ for the
Riemannian gradient, $|\cdot|$ for the length on the tangent space (forgetting the subscript $x$ for simplicity) and
$\|\cdot\|_{L^p}$ for the norm on $ L^p:=L^{p}(M,\mu)$, $1 \leq p\leq +\infty.$  We denote by $Q(x, r)$ the open ball of
center $x\in M $ and radius $r>0$. 
We deal with the Sobolev spaces of order $1$, $W^{1,p}:=W^{1,p}(M)$, where the norm is defined by:
$$ \| f\|_{W^{1,p}(M)} : = \| f\|_p+\|\, |\nabla f|\,\|_{L^p}.$$
We write ${\mathcal S}(M)$ for the Schwartz space on the manifold $M$ and ${\mathcal S}'(M)$ for its dual, corresponding to the set of distributions. Moreover in all this work, ${\bf 1}={\bf 1}_M$ will be used for the constant function, equals to one on the whole manifold.

\subsection{The doubling property}

\begin{df}[Doubling property] Let $M$ be a Riemannian manifold. One says that $M$ satisfies the doubling property $(D)$ if there exists a constant $C_0>0$, such that for all $x\in M,\, r>0 $ we have
\begin{equation*}\tag{$D$}
\mu(B(x,2r))\leq C_0 \mu(B(x,r)).
\end{equation*}
\end{df}

\begin{lem} Let $M$ be a Riemannian manifold satisfying $(D)$ and let $d:=log_{2}C_0$. Then for all $x,\,y\in M$ and $\theta\geq 1$
\begin{equation}\label{eq:d}
\mu(B(x,\theta R))\leq C\theta^{d}\mu(B(x,R))
\end{equation}
There also exists $c$ and $N\geq 0$, so that for all $x,y\in M$ and $r>0$
\be{eq:N} \mu(B(y,r)) \leq c\left(1+\frac{d(x,y)}{r} \right)^N \mu(B(x,r)). \ee
\end{lem} 
\noindent For example, if $M$ is the Euclidean space $M=\R^d$ then $N=0$ and $c=1$. \\
Observe that if $M$ satisfies $(D)$ then
$$ \textrm{diam}(M)<\infty\Leftrightarrow\,\mu(M)<\infty\,\textrm{ (see \cite{ambrosio1})}. $$
Therefore if $M$ is a complete Riemannian manifold satisfying $(D)$ then $\mu(M)=\infty$.

\begin{thm}[Maximal theorem]\label{MIT} (\cite{coifman2})
Let $M$ be a Riemannian manifold satisfying $(D)$. Denote by $\M$ the uncentered Hardy-Littlewood maximal function
over open balls of $M$ defined by
 $$ \M f(x):=\underset{\genfrac{}{}{0pt}{}{B \ \textrm{ball}}{x\in B}} {\sup} \ \frac{1}{\mu(B)}\int_{B} |f| d\mu.  $$
Then for every  $p\in(1,\infty]$, $\M$ is $L^p$-bounded and moreover of weak type $(1,1)$\footnote{ An operator $T$ is of weak type $(p,p)$ if there is $C>0$ such that for any $\alpha>0$, $\mu(\{x;\,|Tf(x)|>\alpha\})\leq \frac{C}{\alpha^p}\|f\|_p^p$.}.
\\
Consequently for $s\in(0,\infty)$, the operator $\M_{s}$ defined by
$$ \M_{s}f(x):=\left[\M(|f|^s)(x) \right]^{1/s} $$
is of weak type $(s,s)$ and $L^p$ bounded for all $p\in(s,\infty]$.
\end{thm}

\subsection{Poincar\'e inequality}
\begin{df}[Poincar\'{e} inequality on $M$] \label{classP} We say that a complete Riemannian manifold $M$ admits a Poincar\'{e} inequality $(P_{q})$ for some $q\in[1,\infty)$ if there exists a constant $C>0$ such that, for every function $f\in W^{1,q}_{loc}(M)$ (the set of compactly supported Lipschitz functions on $M$) and every ball $B$ of $M$ of radius $r>0$, we have
\begin{equation*}\tag{$P_{q}$}
\left(\aver{B}\left|f- \aver{B}f d\mu \right|^{q} d\mu\right)^{1/q} \leq C r \left(\aver{B}|\nabla f|^{q}d\mu\right)^{1/q}.
\end{equation*}
\end{df}
\begin{rem} By density of $C_{0}^{\infty}(M)$ in $W^{1,q}_{loc}(M)$, we can replace $W^{1,q}_{loc}(M)$ by $C_{0}^{\infty}(M)$.
\end{rem}
Let us recall some known facts about Poincar\'{e} inequalities with varying $q$.
 \\
It is known that $(P_{q})$ implies $(P_{p})$ when $p\geq q$ (see \cite{hajlasz4}). Thus, if the set of $q$ such that
$(P_{q})$ holds is not empty, then it is an interval unbounded on the right. A recent result of S. Keith and X. Zhong
(see \cite{KZ}) asserts that this interval is open in $[1,+\infty[$~:

\begin{thm}\label{kz} Let $(M,d,\mu)$ be a doubling and complete Riemannian manifold, admitting a Poincar\'{e} inequality $(P_{q})$, for  some $1< q<\infty$.
Then there exists $\epsilon >0$ such that $(M,d,\mu)$ admits
$(P_{p})$ for every $p>q-\epsilon$.
\end{thm}

\mb We refer the reader to Theorem 5.3.3 of \cite{SC} or Proposition 1.6 of \cite{BB} for the proof of the following consequence.

\begin{prop} \label{prop:poincare}
 Assume that $M$ satisfies $(D)$ and admits a Poincar\'e inequality $(P_q)$ for some $q\in[1,\infty)$. Then there is a constant $c=c(q)$ and $\epsilon>0$ such that for all function $f\in W^{1,q}_{loc}$ 
\be{prop:moy2} \left|f(x)-f(y)\right| \leq c d(x,y) \left[\mathcal{M}_{q-\epsilon}(|\nabla f|)(x) + \mathcal{M}_{q-\epsilon}(|\nabla f|)(y) \right]. \ee
\end{prop}

\subsection{Framework for semigroup of operators} \label{subsec:semigroup}

Let us recall the framework of \cite{DY1, DY}. \\
Let $\omega \in[0,\pi/2)$. We define the closed sector in the complex plane ${\mathbb C}$ by
$$ S_\omega:= \{z\in\C,\ |\textrm{arg}(z)|\leq \omega\} \cup\{0\}$$
and denote the interior of $S_\omega$ by $S_\omega^0$.
We set $H_\infty(S^0_\omega)$ for the set of bounded holomorphic functions $b$ on $S_\omega^0$, equipped with the norm
$$ \|b\|_{H_\infty(S_\omega^0)} := \|b\|_{L^\infty(S_\omega^0)}.$$
Then consider a linear operator $L$. It is said of type $\omega$ if its spectrum $\sigma(L)\subset S_\omega$ and for each $\nu>\omega$, there exists a constant $c_\nu$ such that
$$ \left\|(L-\lambda)^{-1} \right\|_{L^2\to L^2} \leq c_\nu |\lambda|^{-1}$$
for all $\lambda\notin S_\nu$.

\mb We refer the reader to \cite{DY1} and \cite{Mc} for more details concerning holomorphic calculus of such operators. In particular, it is well-known that $L$ generates a holomorphic semigroup $(\A_z:=e^{-zL})_{z\in S_{\pi/2-\omega}}$. Let us detail now some assumptions, we make on the semigroup.

\mb Assume the following conditions: there exist a positive real $m>0$ and $\delta>1$ with \footnote{Usually (see (\cite{DY1, DY}), we just require that $\delta>0$. In our work, we will use Poincar\'e inequality and so we need to compensate for a power $1$ of the distance. That is why, we require $\delta>1$.}
\begin{itemize}
 \item For every $z\in S_{\pi/2-\omega}$, the linear operator $\A_z:=e^{-zL}$ is given by a kernel $a_z$ satisfying
\be{eq:pointwise} \left|a_{z}(x,y)\right|\lesssim \frac{1}{\mu(B(x,|z|^{1/m}))} \left(1+\frac{d(x,y)}{|z|^{1/m}}\right)^{-d-2N-\delta} \ee
where $d$ is the homogeneous dimension of the space (see (\ref{eq:d})) and $N$ is the other dimension parameter (see (\ref{eq:N})); $N\geq0$ could be equal to $0$.
 \item The operator $L$ has a bounded $H_\infty$-calculus on $L^2$. That is, there exists $c_\nu$ such that for $b\in H_\infty(S^0_\nu)$, we can define $b(L)$ as a $L^2$-bounded linear operator and
\be{eq:holocal} \|b(L)\|_{L^2\to L^2} \leq c_\nu \|b\|_\infty. \ee 
\end{itemize}

\begin{rem} \label{rem:holo} The assumed bounded $H_\infty$-calculus on $L^2$ allows us to deduce some extra properties (see \cite{DY} and \cite{Mc})~:
\begin{itemize}
 \item Due to the Cauchy formula for complex differentiation, poinwise estimate (\ref{eq:pointwise}) still holds for the kernel of $(tL)^k e^{-tL}$ with $t>0$ and $k\in {\mathbb N}$.
 \item For any holomorphic function $\psi \in H(S_\nu^0)$ such that for some $s>0$, $ |\psi(z)|\lesssim \frac{|z|^s}{1+|z|^{2s}},$
the quadratic functional 
$$ f \rightarrow \left( \int_0^\infty \left|\psi(tL) f \right|^2 \frac{dt}{t} \right)^{1/2}$$
is $L^2$-bounded.
\end{itemize}
\end{rem}

\mb Moreover we need another assumption on the semigroup, which concerns a square estimate on the gradient of the semigroup. We assume that for every integer $k\geq 0$ the square functional
\be{ass:square} f \rightarrow \left( \int_0^\infty \left|t^{1/m}\nabla (tL)^k e^{-tL}(f) \right|^2 \frac{dt}{t} \right)^{1/2} \ee
is bounded on $L^2$. \\
It is interesting to remark that excepted (\ref{ass:square}), we do not require regularity assumptions on the heat kernel.

\begin{rem} \label{rem:square} We claim that Assumption (\ref{ass:square}) is satisfied under the $L^2$-boundedness of the Riesz transform ${\mathcal R}:= \nabla L^{-1/m}$. \\
Indeed if ${\mathcal R}$ is $L^2$-bounded, then it admits $l^2$-valued estimates, which yields
$$ \left( \int_0^\infty \left| {\mathcal R}  (tL)^{k+1/m} e^{-tL}(f) \right|^2 \frac{dtd\mu }{t} \right)^{1/2} \leq \|{\mathcal R}\|_{L^2\to L^2} \left( \int_0^\infty \left| (tL)^{k+1/m} e^{-tL}(f) \right|^2 \frac{dtd\mu }{t} \right)^{1/2}.$$
This gives the desired result 
$$\left( \int_0^\infty \left| t^{1/m} \nabla (tL)^k e^{-tL}(f) \right|^2 \frac{dtd\mu }{t} \right)^{1/2} \lesssim \|f\|_{L^2},$$
thanks to Remark \ref{rem:holo}.
\end{rem}

\begin{rem} Using Stein's complex interpolation theorem, it is known that a boundedness of 
$$ f \rightarrow \left( \int_0^\infty \left|s^{1/m}\nabla e^{-sL}(f) \right|^2 \frac{dt}{t} \right)^{1/2} $$
on $L^p$ for every $p$ belonging to a neighborhood of $2\in(1,\infty)$ implies the $L^2$-boundedness of (\ref{ass:square}) for $k\in\N$. We refer the reader to Step 7 of the proof of Theorem 6.1 in \cite{A} for a detailed proof of such result.
\end{rem}

\subsection{BMO space} \label{subsec:BMO}

 According to \cite{DY}, it is well-known that BMO spaces related to semigroups are well-defined as a subspace of 
\be{eq:M} {\mathbb M}:=\bigcup_{x_0\in M} \ \bigcup_{\beta\in(0,\delta)} {\mathbb M}_{x_0,\beta} \ee
with ${\mathbb M}_{x_0,\beta}$ the set of functions $f\in L^1_{loc}$ such that
$$ \|f\|_{{\mathbb M}_{x_0,\beta}}:=\int \frac{|f(x)|}{(1+d(x_0,x))^{2N+\beta} \mu(B(x_0,1+d(x_0,x)))} d\mu(x) <\infty,$$
where $N$ is given by (\ref{eq:N}) and $\delta$ by (\ref{eq:pointwise}).

\begin{df}[Definition 2.4 \cite{DY}] \label{def:bmo} A function $f\in {\mathbb M}$ belongs to $BMO_L$ if and only if 
$$ \|f\|_{BMO_L}:= \sup_{t>0} \ \sup_{\genfrac{}{}{0pt}{}{B \textrm{ball}}{r_B^m=t}} \ \frac{1}{\mu(B)}\int_B \left|f-e^{-tL}f\right| d\mu <\infty.$$
 \end{df}

\mb We refer to \cite{DY} and \cite{DY1} for study of precise examples concerning this kind of BMO spaces. We underline that $BMO_{L}$ satisfies to some John-Nirenberg properties and we now state some useful properties:

\begin{prop}[Theorem 2.14 \cite{DY}] \label{prop:Carleson} For any $f\in BMO_L$, the measure
 $$ d\nu(x,t):=\left|t^mL e^{-t^mL} \left(1-e^{-t^mL}\right)\right|^2 \frac{d\mu(x)dt}{t}$$
is a Carleson measure. Moreover for all integer $k\geq 1$ the measure
$$ d\nu_k(x,t):=\left|(t^mL)^k e^{-t^mL} \left(1-e^{-t^mL}\right)\right|^2 \frac{d\mu(x)dt}{t}$$
is a Carleson measure too.
\end{prop}

\dem In \cite{DY}, the proof is only explained for the Laplacian operator $L=-\Delta$ in the Euclidean space with $k=1$. 
However, the authors used technical properties on BMO spaces, which are detailed and proved in the general framework (for example John-Nirenberg inequality). In addition the proof relies on the $L^2$-boundedness of the square function
$$ \left(\int_0^\infty \left|(t^mL)^k e^{-t^mL} \left(1-e^{-t^mL}\right)f \right|^2 \frac{dt}{t}\right)^{1/2},$$
which is a consequence of Assumption (\ref{eq:holocal}), see Remark \ref{rem:holo}. The allowing integer $k\geq 1$ is already appeared in Theorem 9.1 of \cite{HM} in a particular context.
\findem

\begin{rem} In our main result Theorem \ref{thm:principal}, we require the extra assumption 
$$ L({\bf 1}) = L^*({\bf 1}) = 0.$$
Under this assumption, it is known that the classical space BMO (of John-Nirenberg) is included into the new one~:
$$ BMO \hookrightarrow BMO_{L},$$
see Proposition 6.7 \cite{DY1} and Remark 7.6 of \cite{BZ} for a more general study of this question. In addition, this inclusion may be strict, we refer the reader to Proposition 6.8 \cite{DY1} for an example.
\end{rem}

\subsection{Examples of such semigroups} \label{subsec:ex}

In this subsection, we would like to give two examples of situations where all these assumptions are satisfied.

\subsubsection{Second-order elliptic operator}
Let $M=\R^d$ and $A$ be an $d\times d$ matrix-valued function satisfying the ellipticity condition~: there exist two constants $\Lambda\geq \lambda>0$ such that
$$ \forall \xi,\zeta\in \C^d, \qquad \lambda|\xi|^2 \leq Re \left( A\xi \cdot \overline{\xi} \right) \quad  \textrm{and} \quad |A\xi \cdot \overline{\zeta} | \leq \Lambda|\xi||\zeta|.$$ We define the second order divergence form operator
$$ L(f):= -\textrm{div} (A \nabla f).$$
This particular framework was studied for example by P. Auscher in \cite{A},  by S. Hofmann and S. Mayboroda \cite{HM}, ...
We define the second order divergence form operator $L(f)=-div(A \nabla f)$, which can be interpreted in the weak sense via sesquilinear form. Since $L$ is maximal accretive, it admits a bounded $H_\infty$-calculus on $L^2(\R^d)$. Moreover when $A$ has real entries or when the dimension $d\in\{1,2\}$, then the operator $L$ generates an analytic semigroup on $L^2$ with a heat kernel satisfying Gaussian upper-bounds. In this case all our assumptions are verified with $m=2$ (see \cite{A} for estimates on the gradient of the semigroup).

\subsubsection{Laplacian operators on a manifold} \label{subsub:manifold} 
Let $M$ be a doubling connected non-compact Riemannian manifold and consider $L=-\Delta$ the positive Laplace-Beltrami operator. Let $p_t$ the heat kernel of $e^{-tL}$. It is well-known that on-diagonal upper bound
\be{assm} \sup_{x\in M} |p_t(x,x)| \lesssim t^{-d/2} \ee
self-improves into a Gaussian off-diagonal bound (which implies (\ref{eq:pointwise}) with $m=2$). This can be seen in several ways, by using a perturbation method of Davies, or the integrated maximum principle, or the finite propagation speed of solutions to the wave equation (see \cite{C, G, SC, Si}) ... 
Moreover, we know that (\ref{assm}) implies the $L^p$-boundedness of the Riesz transform for every $p\in (1,2]$ (see Theorem 1.1 \cite{CD1}). \\
In addition, by integration by parts we easily obtain the $L^2$-boundedness of the Riesz transform and so
(\ref{ass:square}) is verified, thanks to Remark \ref{rem:square}.
We refer the reader to \cite{ACDH} for more details concerning these kind of assumptions and how they are related between them. \\
Concerning the extension to a complex semigroup, we refer the reader to \cite{Da} (Lemma 2). Since $L$ is a non-negative self-adjoint operator, the semigroup $e^{-zL}$ is holomorphic in $S_{\pi/2}$ and $L$ admits a $H_\infty$-bounded holomorphic calculus.

\gb Consequently, for example if the manifold $M$ has nonnegative Ricci curvature all our assumptions are satisfied.

\section{A T(1) theorem for semigroup}

\subsection{Assumptions and statement}

Before stating our $T(1)$ theorem, we have to assume some properties on a generic operator $T$ (in order to replace the usual Calder\'on-Zygmund properties on the kernel in the classical T(1) theorem). Let $T$ be a continuous operator acting from ${\mathcal S}(M)$ into ${\mathcal S}'(M)$.
We assume that $T$ and $T^*$ satisfies some $L^2-L^2$ off diagonal decay as follows~: there exists an integer $\K\geq 1$ such that for every $s>0$, every ball $Q_1,Q_2$ of radius $r:=s^{1/m}$ and function $f\in L^2(Q_1)$ 
\begin{itemize}
\item if $d(Q_1,Q_2)\geq 2 r$, then we have 
\be{assum:off} \left\|(sL)^\K e^{-sL} T(f)\right\|_{L^2(Q_2)} \lesssim \left(1+\frac{d(Q_1,Q_2)}{r}\right)^{-d-2N-\delta} \|f\|_{L^2(Q_1)} \ee
and the dual estimates
\be{assum:off-dual} \left\|(sL^*)^\K e^{-sL^*} T^*(f)\right\|_{L^2(Q_2)} \lesssim \left(1+\frac{d(Q_1,Q_2)}{r}\right)^{-d-2N-\delta} \|f\|_{L^2(Q_1)}. \ee
\item if $d(Q_1,Q_2)\leq 2 r$, then we have
\be{assum:weaklybounded} \left\|(sL)^{\K} e^{-sL} T(e^{-sL} f)\right\|_{L^2(Q_2)} + \left\|(sL^*)^{\K}e^{-sL^*} T(e^{-sL^*} f)\right\|_{L^2(Q_2)} 
\lesssim \|f\|_{L^2(Q_1)}. \ee
We call this property the {\it weak boundedness} of $T$.
\end{itemize}

\begin{rem} \label{rem:weakly} We claim that (\ref{assum:weaklybounded}) self-improves in the following one~: for all integer $k\geq 0$ (with $d(Q_1,Q_2)\leq 2r$)
$$ \left\|(sL)^{\K} e^{-sL} T((sL)^k e^{-sL} f)\right\|_{L^2(Q_2)} + \left\|(sL^*)^{\K}e^{-sL^*} T((sL)^k e^{-sL} f)\right\|_{L^2(Q_2)} \lesssim \|f\|_{L^2(Q_1)}. $$
Let us explain how can we obtain this improvement for $T$. We chose a bounded covering $(R_l)_l$ of the manifold  by balls of radius $r:=(s/2)^{1/m}$. Applying (\ref{assum:weaklybounded}) and (\ref{assum:off}) to the function $(sL)^k e^{-sL/2}(f)$, we get
$$ \left\|(sL)^{\K} e^{-sL} T((sL)^k e^{-sL} f)\right\|_{L^2(Q_2)} \lesssim \sum_{l} \left(1+\frac{d(Q_2,R_l)}{r}\right)^{-d-2N-\delta} \|(sL)^k e^{-sL/2}(f)\|_{L^2(R_l)}.$$
Then, due to the off-diagonal decay of the derivative of the semigroup, it comes
\begin{align*} 
\left\|(sL)^{\K} e^{-sL} T((sL)^k e^{-sL} f)\right\|_{L^2(Q_2)} & \\
& \hspace{-3cm} \lesssim \sum_{l} \left(1+\frac{d(Q_2,R_l)}{r}\right)^{-d-2N-\delta} \left(1+\frac{d(Q_1,R_l)}{r}\right)^{-d-2N-\delta} \|f\|_{L^2(Q_1)} \\
& \hspace{-3cm} \lesssim \|f\|_{L^2(Q_1)}.
\end{align*}
\end{rem}

\begin{rem} The {\it weak boundedness} property of the operator $T$ is obviously a necessary condition for the $L^2$-boundedness.
\end{rem}

\begin{thm} \label{thm:principal} Assume that the Riemannian manifold $M$ satisfies Poincar\'e inequality $(P_2)$ and has a infinite measure $\mu(M)=\infty$. Suppose the existence of an operator $L$ such that $L({\bf 1})=0=L^*({\bf 1})$ and such that the corresponding semigroup verifies the assumptions of Subsection \ref{subsec:semigroup}. Let $T$ be a continuous operator weakly continuous from $L^2$ into $L^2$ satisfying (\ref{assum:off}), (\ref{assum:off-dual}) and (\ref{assum:weaklybounded}). 
\begin{itemize}
 \item If $T({\bf 1})\in BMO_L$ and $T^*({\bf 1})\in BMO_{L^*}$ then $T$ admits a $L^2$-bounded extension. Moreover $\|T\|_{L^2\to L^2}$ is only controlled by the implicit constants and $\|T({\bf 1})\|_{BMO_L}+\|T^*({\bf 1})\|_{BMO_{L^*}}$ (and not by the weak continuity from $L^2$ to $L^2$). 
 \item Let assume $\K=1$ in (\ref{assum:off}), (\ref{assum:off-dual}) and (\ref{assum:weaklybounded}). If $T$ admits a continuous extension on $L^2$ then $T({\bf 1})$ belong to $BMO_L$ and $T^*({\bf 1})$ to $BMO_{L^*}$.
\end{itemize}
\end{thm}

\begin{rem} In the previous statement, we implicitly assume that $T({\bf 1})$ and $T^*({\bf 1})$ are well-defined and belong to ${\mathbb M}$. We do not deal with this specific problem as in our applications, we can work with some kind of ``truncations'' $T_{\epsilon}$ of $T$ where $T_\epsilon({\bf 1})$ and $T_\epsilon^*({\bf 1})$ will be well-defined.
\end{rem}

\begin{rem} For the second part of the theorem, for another integer $\K\geq 1$, if $T$ admits a continuous extension on $L^2$ then $T({\bf 1})$ and $T^*({\bf 1})$ belong to another kind of $BMO_L$ space (which requires more cancellation) and satisfies
$$ \sup_{t>0} \ \sup_{\genfrac{}{}{0pt}{}{Q \textrm{ball}}{r_Q^m=t}} \ \frac{1}{\mu(Q)}\int_Q \left|(1-e^{-tL})^\K T({\bf 1})\right| d\mu <\infty $$
and
$$ \sup_{t>0} \ \sup_{\genfrac{}{}{0pt}{}{Q \textrm{ball}}{r_Q^m=t}} \ \frac{1}{\mu(Q)}\int_Q \left|(1-e^{-tL^*})^\K T^*({\bf 1})\right| d\mu <\infty. $$
We let to the reader to check this claim by adapting the proof made for the second part of Theorem 3.2.
\end{rem}

\begin{cor} \label{cor} Under the assumptions of Theorem \ref{thm:principal}, for all exponent $p\in (1,\infty)$ there is a constant $c_p$ such that
$$ \forall f\in {\mathcal S}, \qquad \|T(f)\|_{L^p} \leq c_p \|f\|_{L^p}.$$
\end{cor}

\dem From Theorem \ref{thm:principal}, we know that the operator $T$ admits a continuous extension in $L^2$. To deduce an $L^p$-boundedness from an $L^2$-boundedness, we use the theory of Hardy spaces. \\
More precisely, we refer to \cite{BZ} for an abstract theory of Hardy spaces, which we are going to apply. We have to build an adapted Hardy space $H^1_L$ (step 1), then to prove that our operator $T$ is bounded from $H^1_L$ to $L^1$ (step 2) and then to interpolate the Hardy space with $L^2$ (step 3). \\
{\bf Step 1:} Construction of an adapted Hardy space. \\
For $Q$ a ball of radius $r=s^{1/m}$, we define the operator 
$$ A_Q:= 1-(1-e^{-sL})^\K.$$
These operators are uniformly bounded on $L^2$ (since it corresponds to a finite sum of semigroups).
We refer the reader to \cite{BZ} for an abstract construction of Hardy space $H^1_{L}$, based on a collection of operators $(A_Q)_Q$, indexed by the balls. We refer to \cite{DY1} for more specific results concerning the particular case where $A_Q$ is defined by a semigroup. Indeed, the Hardy space is defined by atomic decomposition. By definition for a ball $Q$, an atom (relatively to the ball $Q$) is a function $m=f-A_Q(f)$ where $f\in L^2(Q)$ is $L^1$-normalized in $L^2$.

\mb
{\bf Step 2:} Boundedness of $T$ from $H^1_L$ to $L^1$. \\
To check that $T$ is bounded on the implicit Hardy space $H^1_L$, it just suffices to prove some $L^2-L^2$ off-diagonal decay (see Theorem 4.2 of \cite{BZ})~: for every ball $Q$ or radius $r=s^{1/m}$ and integer $i\geq 1$
\be{hardyamontrer}
\left(\frac{1}{\mu(2^iQ)} \int_{2^iQ\setminus 2^{i-1}Q} \left|T\left(1-e^{-sL}\right)^\K(f)\right|^2 d\mu \right)^{1/2} \lesssim \gamma(i) \|f\|_{L^2(Q)} \ee
with fast decreasing coefficients $\gamma(i)$ such that
$$\sum_{i} 2^{id}\gamma(i)<\infty$$ and
\be{hardyamontrer2}
\left(\frac{1}{\mu(2Q)} \int_{2Q} \left|T\left(1-e^{-sL}\right)^\K(f)\right|^2 d\mu \right)^{1/2} \lesssim  \|f\|_{L^2(Q)}. \ee
These off-diagonal estimates are a consequence of (\ref{assum:off-dual}). Indeed off-diagonal decay (\ref{assum:off-dual}) can be improved as follows: for every $s>0$, every ball $Q_1,Q_2$ of radius $r:=s^{1/m}$ and function $f\in L^2(Q_1)$ then we have
\be{assum:off2} \left\|(sL^*)^\K e^{-sL^*} T^*(f)\right\|_{L^2(Q_2)} \lesssim \left(1+ \frac{d(Q_1,Q_2)}{r}\right)^{-d-2N-\delta} \|f\|_{L^2(Q_1)}. \ee
Actually if $d(Q_1,Q_2)\geq 2 r$, this is (\ref{assum:off-dual}) and if $d(Q_1,Q_2)\leq 2 r$ then we use $L^2$-boundedness of $T^*$ and $(sL^*)^\K e^{-sL^*}$. \\
From (\ref{assum:off2}), we deduce by duality that 
\be{assum:off3} \left\|T(sL)^\K e^{-sL}(f)\right\|_{L^2(Q_2)} \lesssim \left(1+ \frac{d(Q_1,Q_2)}{r}\right)^{-d-2N-\delta} \|f\|_{L^2(Q_1)}. \ee
To deduce (\ref{hardyamontrer}) when $d(Q_1,Q_2)\geq 2r$, we differentiate the semigroup as follows
\begin{align*}
\left(\frac{1}{\mu(Q_2)} \int_{Q_2} \left|T\left(1-e^{-sL}\right)^{\K}(f)\right|^2 d\mu \right)^{1/2} & \\
& \hspace{-4cm} \lesssim \left(\frac{1}{\mu(Q_2)} \int_{Q_2} \left|\int_0^s \cdots \int_0^s  TL^\K e^{-(u_1+\cdots +u_\K)L}(f) du_1\cdots du_\K\right|^2 d\mu \right)^{1/2}.
\end{align*}
For each $u_i\in(0,s]$, we have to consider $(B_k)_k$ and $(B_l)_l$ a bounded covering of $Q_2$ and $Q_1$ by balls of radius $r_u=(u_1+\cdots + u_\K)^{1/m}$ and then to apply (\ref{assum:off3}):
\begin{align*}
\int_{Q_2} \left| r_u^{m\K} TL^\K e^{-(u_1+\cdots +u_\K)L}(f) \right|^2 d\mu & \leq \sum_k \int_{B_k} \left| r_u^{m\K} TL^\K e^{-(u_1+\cdots +u_\K)L}(f)\right|^2 d\mu \\
& \hspace{-4cm} \lesssim \sum_{k} \left( \sum_l \left(1+ \frac{d(B_k,B_l)}{r_u}\right)^{-d-2N-\delta} \left(\frac{\mu(B_k)}{\mu(B_l)}\right)^{1/2}\left(\int_{B_l} |f|^2 d\mu\right)^{1/2}\right)^2 \\
& \hspace{-4cm} \lesssim \sum_{k} \left(\sum_l \left(1+ \frac{d(B_k,B_l)}{r_u}\right)^{-d-N-\delta} \left(\int_{B_l} |f|^2 d\mu\right)^{1/2}\right)^2 \\
& \hspace{-4cm} \lesssim \left(\int_{Q_1} |f|^2 d\mu\right) \sum_{k,l} \left(1+ \frac{d(B_k,B_l)}{r_u}\right)^{-2d-2N-2\delta},
\end{align*}
where we have used the doubling property and (\ref{eq:N}). Then we note that $d(B_k,B_l)\geq d(Q_1,Q_2)$ and due to the doubling property,
\begin{align}
 \sum_{k} 1 & \lesssim \left(\frac{r}{r_u}\right)^d \sum_{k} \left(\frac{r}{r_u}\right)^{-d} \lesssim \left(\frac{r}{r_u}\right)^d \sum_{k} \frac{\mu(B_k)}{\mu(\frac{r}{r_u}B_k)} \nonumber \\
&  \lesssim \left(\frac{r}{r_u}\right)^d \sum_{k} \frac{\mu(B_k)}{\mu(Q_2)} \nonumber \\
&  \lesssim \left(\frac{r}{r_u}\right)^d. \label{eq:argument}
\end{align}
So it follows with a similar reasoning for the sum over $l$ that
\begin{align*}
\int_{Q_2} \left| r_u^{m\K} TL^\K e^{-(u_1+\cdots +u_\K)L}(f) \right|^2 d\mu & \lesssim \left(\int_{Q_1} |f|^2 d\mu\right) \left(1+ \frac{d(Q_1,Q_2)}{r_u}\right)^{-2d-2N-2\delta} \left(\frac{r}{r_u}\right)^{2d}.
\end{align*}
We conclude that
\begin{align*}
 \left(\frac{1}{\mu(Q_2)} \int_{Q_2} \left| r_u^{m\K} TL^\K e^{-(u_1+\cdots +u_\K)L}(f)\right|^2 d\mu\right)^{1/2} & \\
& \hspace{-5cm} \lesssim \left(1+\frac{d(Q_1,Q_2)}{r_u}\right)^{-d-N-\delta} \left(\frac{r}{r_u}\right)^d \left(\frac{\mu(Q_1)}{\mu(Q_2)}\right)^{1/2} \left(\frac{1}{\mu(Q_1)} \int_{Q_1} |f|^2 d\mu\right)^{1/2}.
\end{align*}
Hence, for $d(Q_1,Q_2)\geq 2r$
\begin{align} 
\left(\frac{1}{\mu(Q_2)} \int_{Q_2} \left|T\left(1-e^{-sL}\right)(f)\right|^2 d\mu \right)^{1/2} & \\
& \hspace{-5cm} \lesssim \left(\frac{\mu(Q_1)}{\mu(Q_2)}\right)^{1/2} \left(\frac{1}{\mu(Q_1)} \int_{Q_1} |f|^2 d\mu\right)^{1/2} \int_{0}^{\K s} \left(1+\frac{d(Q_1,Q_2)}{u^{1/m}}\right)^{-d-N-\delta} \left(\frac{r}{u^{1/m}}\right)^d \frac{du}{u}  \nonumber \\
& \hspace{-5cm} \lesssim  \left(\frac{\mu(Q_1)}{\mu(Q_2)}\right)^{1/2} \left(1+\frac{d(Q_1,Q_2)}{r}\right)^{-d-N-\delta}  \left(\frac{1}{\mu(Q_1)} \int_{Q_1} |f|^2 d\mu\right)^{1/2} \nonumber \\
& \hspace{-5cm} \lesssim \left(1+\frac{d(Q_1,Q_2)}{r}\right)^{-d-\delta}  \left(\frac{1}{\mu(Q_1)} \int_{Q_1} |f|^2 d\mu\right)^{1/2}. \label{eq:sollu}
\end{align}
We have used at the first line that
$$ \left| \left\{(u_1,\cdots,u_\K),\ u_1+ \cdots + u_\K=u\right\} \right| \lesssim u^{\K-1}.$$
If $d(Q_1,Q_2)\leq 2r$, (\ref{eq:sollu}) still holds, indeed we have not to differentiate and just invoke the $L^2$-boundedness of the different appearing operators. \\
So we have proved that for every balls $Q_1,Q_2$ of radius $r$, (\ref{eq:sollu}) holds. Now we deduce (\ref{hardyamontrer}) as follows: we consider $Q$ a ball of radius $r$ and $(\tilde{Q}_k)_k$ a bounded covering of $2^iQ\setminus 2^{i-1} Q$ by balls of radius $r$ and then apply (\ref{eq:sollu}):
\begin{align*}
\left(\frac{1}{\mu(2^iQ)} \int_{2^iQ\setminus 2^{i-1}Q_1} \left|T\left(1-e^{-sL}\right)(f)\right|^2 d\mu\right)^{1/2} & \\
& \hspace{-3cm} \leq  \left(\frac{1}{\mu(2^iQ)} \sum_k \int_{\tilde{Q}_k } \left| T\left(1-e^{-sL}\right)(f)\right|^2 d\mu\right)^{1/2} \\
&  \hspace{-3cm} \lesssim  \left( \sum_k \left(1+\frac{d(\tilde{Q}_k,Q)}{r}\right)^{-2d-2\delta} \frac{\mu(\tilde{Q}_k)}{\mu(2^iQ)\mu(Q)}\int_{Q} |f|^2 d\mu\right)^{1/2} \\
&  \hspace{-3cm} \lesssim 2^{-(d+\delta)i} \left(\frac{1}{\mu(Q)}\int_{Q} |f|^2 d\mu\right)^{1/2}.
\end{align*}
We have also obtained (\ref{hardyamontrer}) with coefficients $\gamma(j)$ satisfying
$$ \gamma(i) \lesssim 2^{-(d+\delta)i}.$$
Then Theorem 4.2 of \cite{BZ} yields that $T$ is bounded from the finite atomic Hardy space $H^1_L$.

\mb
{\bf Step 3:} Interpolation between $H^1_L$ and $L^2$. \\
To obtain interpolation results, we apply Theorem 5.3 of \cite{BZ}. Aiming that, off-diagonal estimates of the heat kernel (\ref{eq:pointwise}) imply that
\begin{align*} 
M_\infty(f)(x)& :=\sup_{Q \textrm{ball}} \|A_Q^*(f)\|_{L^\infty(Q)} \\
& \lesssim \sup_{r>0} \int \frac{1}{\mu(B(y,r))} \left(1+\frac{d(x,y)}{r}\right)^{-d-2N-\delta} |f(y)| d\mu(y) \\
& \lesssim \sup_{r>0} \frac{1}{\mu(B(x,r))} \int \left(1+\frac{d(x,y)}{r}\right)^{-d-N-\delta} |f(y)| d\mu(y) \\
& \lesssim {\mathcal M}(f)(x),
\end{align*}
where we have used (\ref{eq:N}) and the fact that $A_Q^*$ can be expanded as a finite sum of semigroups (with a scale equivalent to $r$ the radius of the ball $Q$). \\
Consequently, we know from Theorem 5.3 \cite{BZ} that we can interpolate the Hardy space $H^1_{L}$ associated to the operators $A_Q$ with $L^2$ and regain the intermediate Lebesgue spaces, as intermediate spaces (see \cite{B} for more details concerning a real interpolation result). This concludes the proof of the theorem for $p\in (1,2]$. We get the result for $p\in(2,\infty)$ by duality, since $T$ and $T^*$ satisfy the same assumptions.
 \findem

\mb In addition, we can obtain weighted estimates. Let us recall the definition of Muckenhoupt's weights and Reverse H\"older classes~:

\begin{df} \label{def:poids} A nonnegative function $\omega$ on $M$ belongs to the class ${\mathbb A}_p$ for $p\in (1,\infty)$ if
$$\sup_{Q \textrm{ ball}} \left(\frac{1}{\mu(Q)}\int_Q w d\mu \right) \left( \frac{1}{\mu(Q)}\int_Q \omega^{-1/(p-1)} d\mu \right)^{p-1} <\infty.$$
A nonnegative function $\omega$ on $M$ belongs to the class $RH_q$ for $q\in[1,\infty)$, if there is a constant $C$ such that for every ball $Q\subset X$
$$ \left( \frac{1}{\mu(Q)}\int_Q \omega^q d\mu \right)^{1/q} \leq C \left(\frac{1}{\mu(Q)}\int_Q \omega d\mu \right).$$
Remark that $RH_1$ is the class of all the weights. \\
For $\omega$ a weight and $p\in[1,\infty]$ an exponent, we write $L^p(\omega)$ for the weighted Lebesgue space (associated to the measure $\omega d\mu$).
\end{df}

\mb Following previous Corollary \ref{cor} and Theorem 6.4 of \cite{BZ}, we get the next result.

\begin{cor} \label{cor2} Under the assumptions of Theorem \ref{thm:principal}, for all exponent $p\in (1,2)$, let $\omega$ be a weight so that $\omega \in {\mathbb A}_{p} \cap RH_{(\frac{2}{p})'}$. There is a constant $c_p$ such that
$$ \forall f\in {\mathcal S}(M), \qquad \|T^*(f)\|_{L^p(\omega)} + \|T(f)\|_{L^p(\omega)} \leq c_p \|f\|_{L^p(\omega)}. $$
For all exponent $p\in (2,\infty)$, let $\omega$ be a weight so that $\omega \in {\mathbb A}_{p/2}$. There is a constant $c_p$ such that
$$ \forall f\in {\mathcal S}(M), \qquad \|T^*(f)\|_{L^p(\omega)} + \|T(f)\|_{L^p(\omega)} \leq c_p \|f\|_{L^p(\omega)}.$$ 
\end{cor}

\mb Concerning the condition $\omega \in {\mathbb A}_{p/2}= {\mathbb A}_{p/2}\cap RH_{1}$, we recall (Lemma 4.4 of \cite{AM}) that for $p\in (2,\infty)$
\be{poidseq}
\omega \in {\mathbb A}_{p/2} \cap RH_{1} \Longleftrightarrow \omega^{1-p'}\in {\mathbb A}_{p'}\cap RH_{(\frac{2}{p'})'}
\ee
and the duality for the weighted Lebesgue spaces goes as follows: for $T$ an operator, we recall the fact that $T^*$ is the adjoint of $T$ related to the measure $\mu$ ; so for $p\in(1,\infty)$
$$ \textrm{ $T$ is $L^p(\omega)$-bounded} \ \Longleftrightarrow \ \textrm{ $T^*$ is $L^{p'}(\omega^{1-p'})$-bounded}.$$

\subsection{proof of Theorem \ref{thm:principal}}

We devote this subsection to the proof of our main result, Theorem \ref{thm:principal}.
Aiming that, let us define a function : let $\phi$ be the following function defined on $\R^+$ by:
\begin{align*}
\phi(x):= -\int_x^\infty t^\K e^{-t} (1-e^{-t})^2 e^{-t} dt.
\end{align*}
Then $\phi$ is a finite sum of polynomial term, multiplied by decreasing exponentials. So it can be extended in $\C$ and becomes a bounded and holomorphic function in any sector $S_\omega$ with $\omega<\pi/2$ and we have
\be{eq:phi} \phi'(x)=x^\K e^{-x}(1-e^{-x})^2 e^{-x}.\ee

\begin{lem} \label{lem:limit} For $f\in L^2(X)$, we have the following strong convergence in $L^2$~:
\begin{itemize}
 \item  $e^{-tL}(f) \xrightarrow[t\to 0]{} f$ 
 \item  for an integer $k\geq 1$, $(tL)^k e^{-tL}(f) \xrightarrow[t\to 0]{} 0$
 \item for an integer $k$, $(tL)^k e^{-tL}(f) \xrightarrow[t\to \infty]{} 0$.
\end{itemize}
\end{lem}

\dem The two first points are due to the $L^2$-continuity of the semigroup. For the third point with $k=0$, we remark that the pointwise bound (\ref{eq:pointwise}) yields
$$ \sup_{t>1} \left|e^{-tL}(f)\right| \lesssim {\mathcal M}(f) \in L^2.$$
As for every $x\in M$ and $t\geq 1$
\begin{align*}
 \left|e^{-tL}(f)(x)\right|& \leq \frac{1}{\mu(B(x,t^{1/m}))} \int \left(1+\frac{d(x,y)}{t^{1/m}}\right)^{-n-2N-\delta} |f(y)| d\mu(y) \\
 & \lesssim \frac{1}{\mu(B(x,t^{1/m}))^{1/2}} \|f\|_{L^2} \xrightarrow[t\to \infty]{} 0.
\end{align*}
At the second inequality, we have used the homogeneous type of the manifold $M$ and then at the last inequality the infinite measure of $M$. As a consequence, we get that $e^{-tL}(f)$ strongly converges to $0$ in $L^2$ when $t\to \infty$. Then for $k\geq 1$, the limit follows from the uniform $L^2$ boundedness of $(tL)^k e^{-tL}$.
\findem

\begin{prop} \label{prop:sobolev} We have a Sobolev inequality, relatively to $L$~: for a large enough integer $M$, for all $t>0$ and ball $Q$ of radius $r=t^{1/m}$
$$ \|f\|_{L^\infty(Q)} \lesssim \inf_{Q} \ {\mathcal M}_2\left[ \left(1+tL\right)^M f \right].$$
More precisely, for some $\delta>1$ (indeed $\delta$ is the one introduced in (\ref{eq:pointwise}))
$$\|f\|_{L^\infty(Q)} \lesssim  \sum_{i\geq 0} 2^{-i\delta} \left( \frac{1}{\mu(2^iQ)} \int_{2^iQ} \left|\left(1+tL\right)^M f\right|^2 d\mu \right)^{1/2}.$$
\end{prop}

\dem It suffices to prove that
$$ \|\left(1+tL\right)^{-M} f\|_{L^\infty(Q)} \lesssim \inf_{Q} {\mathcal M}_2(f)$$
(and similarly for the second estimate), which will be provided as soons as we will prove the pointwise inequality
\be{eq:sobamontrer} \left|r_t(x,y)\right| \leq \frac{1}{\mu(B(x,r))} \left(1+\frac{|x-y|}{r}\right)^{-d-2N-\delta} \ee
where $r_t$ is the kernel of the resolvant $\left(1+tL\right)^{-M}$ and $N,d$ are introduced in (\ref{eq:d}) and (\ref{eq:N}). \\
Aiming that, we decompose the resolvant with the semigroup as follows (up to a numerical constant)
$$ \left(1+tL\right)^{-M} = \int_{0}^\infty s^Me^{-s(1+tL)} \frac{ds}{s}.$$
Using the estimates of the heat kernel (\ref{eq:pointwise}), we obtain
$$ \left|r_t(x,y)\right| \lesssim  \int_{0}^\infty s^Me^{-s} \frac{1}{\mu(B(x,rs^{1/m}))} \left(1+\frac{d(x,y)}{rs^{1/m}}\right)^{-d-2N-\delta} \frac{ds}{s}.$$
The integral for $s\in(0,1]$ is bounded by
\begin{align*}
 \int_{0}^1 s^Me^{-s} \frac{1}{\mu(B(x,rs^{1/m}))} \left(1+\frac{d(x,y)}{rs^{1/m}}\right)^{-d-2N-\delta} \frac{ds}{s} & \\ 
& \hspace{-4cm} \lesssim \int_{0}^1 s^M \frac{1}{\mu(B(x,rs^{1/m}))} \left(1+\frac{d(x,y)}{r}\right)^{-d-2N-\delta} \frac{ds}{s} \\
& \hspace{-4cm} \lesssim \left(1+\frac{d(x,y)}{r}\right)^{-d-2N-\delta} \int_{0}^1 s^{M-1} \frac{1}{\mu(B(x,rs^{1/m}))} ds.
\end{align*}
Using the doubling property of the measure, it follows that
$$ \mu(B(x,r)) \lesssim \mu(B(x,rs^{1/m})) s^{-d/m}.$$
So for $M$ large enough, it comes
$$  \int_{0}^1 s^Me^{-s} \frac{1}{\mu(B(x,rs^{1/m}))} \left(1+\frac{d(x,y)}{rs^{1/m}}\right)^{-d-2N-\delta} \frac{ds}{s} \lesssim \left(1+\frac{d(x,y)}{r}\right)^{-d-2N-\delta} \frac{1}{\mu(B(x,r))},$$
which yields the desired inequality for this first term. \\
Concerning the second part, with $M> (d+2N+\delta)/m$, we have 
\begin{align*}
 \int_{1}^\infty s^M e^{-s} \frac{1}{\mu(B(x,rs^{1/m}))} \left(1+\frac{d(x,y)}{rs^{1/m}}\right)^{-d-2N-\delta} \frac{ds}{s} & \\
& \hspace{-4cm} \lesssim \int_{1}^\infty s^{-M} \frac{1}{\mu(B(x,rs^{1/m}))} \left(1+\frac{d(x,y)}{rs^{1/m}}\right)^{-d-2N-\delta} \frac{ds}{s} \\
& \hspace{-4cm} \lesssim \frac{1}{\mu(B(x,r))} \int_{1}^\infty s^{-M}  \left(1+\frac{d(x,y)}{rs^{1/m}}\right)^{-d-2N-\delta} \frac{ds}{s} \\
& \hspace{-4cm} \lesssim \frac{1}{\mu(B(x,r))}  \left(1+\frac{d(x,y)}{r} \right)^{-d-2N-\delta},
\end{align*}
since for every $v>0$
$$ \int_{1}^\infty s^{-M}\left(1+\frac{v}{s^{1/m}}\right)^{-d-2N-\delta} \frac{ds}{s} \simeq \left(1+v\right)^{-d-2N-\delta},$$
by dividing the integral for $s\leq v^m$ and $s\geq v^m$. This concludes the proof of (\ref{eq:sobamontrer}). Then
we have
\begin{align*}
 \|\left(1+tL\right)^{-M} f\|_{L^\infty(Q)} & \leq  \sup_{x\in Q} \ \int \left|r_t(x,y) \right| \left|f(y)\right| d\mu(y) \\
& \leq \sup_{x \in Q} \  \frac{1}{\mu(B(x,r))} \int \left(1+\frac{|x-y|}{r}\right)^{-d-2N-\delta} \left|f(y)\right| d\mu(y) \\
& \lesssim \sup_{x \in Q} \ \frac{1}{\mu(Q)} \int \left(1+\frac{d(y,Q)}{r}\right)^{-d-2N-\delta} \left|f(y)\right| d\mu(y)\\
& \lesssim \sum_{i\geq 0}  \frac{1}{\mu(Q)} 2^{-(d+2N+\delta)i} \mu(2^iQ)^{1/2} \left\|f\right\|_{L^2(2^iQ)} \\
& \lesssim \sum_{i\geq 0}  2^{-\delta i} \left(\inf_{Q} {\mathcal M}_2 (f) \right),
\end{align*}
where we have used the doubling property (\ref{eq:d}). This completes the proof of the proposition.
\findem

\gb
{\bf Proof of the first part of Theorem \ref{thm:principal}}. \\
The proof is quite long, so we divide it in two steps. First as we use duality, we precise that Assumptions of Subsection \ref{subsec:semigroup} assumed for $L$ still hold for the adjoint $L^*$ (as they are invariant by duality).

\mb
{\bf First step:} Reduction to well-localized operators. \\
Let $f\in L^2$ compactly supported, then in weak $L^2$-sense we have
$$ \phi(0)^2 T(f) = \lim_{s\to 0} \phi(sL) T \phi(sL) (f) $$
and
$$ 0 = \lim_{s\to \infty} \phi(sL) T \phi(sL) (f),$$
thanks to Lemma \ref{lem:limit}. We can also write
$$ T(f) =\lim_{\genfrac{}{}{0pt}{}{\epsilon\to 0}{R\to \infty}} \left[\phi(sL)T\phi(sL) (f)\right]^{R}_\epsilon,$$
which gives
$$ T(f) = \lim_{\genfrac{}{}{0pt}{}{\epsilon\to 0}{R\to \infty}} \int_{\epsilon}^R \left( \left[s\frac{d}{ds} \phi(sL)\right]T\phi(sL) (f) + \phi(sL) T \left[s\frac{d}{ds} \phi(sL)\right] f \right)\frac{ds}{s}.$$
Since the definition of $\phi$, it comes
$$ s\frac{d}{ds} \phi(sL) = (sL)^\K e^{-sL}(1-e^{-sL})^2 sLe^{-sL}$$
in order that
\begin{align*}
 T(f) = & \lim_{\genfrac{}{}{0pt}{}{\epsilon\to 0}{R\to \infty}} \int_{\epsilon}^R \left[ (sL)^\K e^{-sL}(1-e^{-sL})^2sLe^{-sL} T\phi(sL) (f) \right. \\
& \hspace{1cm} \left. + \phi(sL) T (sL)^\K e^{-sL}(1-e^{-sL})^2sLe^{-sL} f \right] \frac{ds}{s}.
\end{align*}
The second term can be seen as the adjoint of the first term by changing $L$ with $L^*$ (and $T$ with $T^*$). As our assumptions are invariant by duality, it is also sufficient to prove a uniform bound with respect to $\epsilon$ and $R$ of $\|U_{\epsilon,R}(f)\|_{L^2}$ with
$$ U_{\epsilon,R}(f):= \int_{\epsilon}^R (1-e^{-sL})sL e^{-sL} T_s \phi(sL) (f) \frac{ds}{s},$$
where we set $T_s:=(sL)^\K e^{-sL}(1-e^{-sL}) T$. \\
According to the assumptions, it is easy to check that $T_s$ satisfies the two following properties:
\begin{itemize}
 \item The map $(s,x)\to T_{s^m}({\bf 1})(x)$ defines a Carleson measure.\\
This comes from $T_s({\bf 1}) = (sL)^\K e^{-sL}(1-e^{-sL}) T({\bf 1}) $ and Proposition \ref{prop:Carleson}.
 \item For all $s$, the operator $T_s$ has off-diagonal decays at the scale $s$~: for every $L^2$ function $f,g$ supported on ball $Q_1,Q_2$ of radius $r:=s^{1/m}$, we have~:
\be{eq:prop2}  \left| \langle T_s(e^{-sL} f),g \rangle \right| \lesssim \left(1+\frac{d(Q_1,Q_2)}{r}\right)^{-d-2N-\delta} \|f\|_{L^2(Q_1)} \|g\|_{L^2(Q_2)}. \ee
This comes from the observation: $\langle T_s(\phi(sL)f),g \rangle = \langle (sL)^\K (1-e^{-sL}) e^{-sL} T(\phi(sL)f),g\rangle$ and Assumption (\ref{assum:off}) together with arguments of Remark \ref{rem:weakly} (if $d(Q_1,Q_2)\geq 2 r$) and Assumption (\ref{assum:weaklybounded}) (if $d(Q_1,Q_2)\leq 2 r$).
\end{itemize}

\mb
{\bf Second step :} Study of $U_{\epsilon,R}$.\\
We decompose $U_{\epsilon,R}$ as follows
\begin{align}
 U_{\epsilon,R}(f) & = \int_{\epsilon}^R (1-e^{-sL})Le^{-sL} \left[T_s({\bf 1}) \phi(sL) (f)\right] ds \nonumber \\
 & \ \ + \int_{\epsilon}^R (1-e^{-sL})Le^{-sL} \left[T_s \phi(sL) - T_s({\bf 1}) \phi(sL)\right](f) ds \nonumber \\
 & := \textrm{Main}(f) + \textrm{Error}(f). \label{eq:error}
\end{align}
The main term is  controled by duality: in pairing with any $L^2$-function $g$
\begin{align*}
 \left| \langle \textrm{Main}(f),g \rangle \right| & = \left| \int_{\epsilon}^R \langle (1-e^{-sL})sLe^{-sL} \left[T_s({\bf 1}) \phi(sL) (f),g \rangle \right] \frac{ds}{s} \right| \\
& = \left| \int_{\epsilon}^R \langle T_s({\bf 1}) \phi(sL) (f),(1-e^{-sL^*})sL^*e^{-sL^*} g \rangle \frac{ds}{s} \right| \\
& \leq \left\| T_s({\bf 1}) \phi(sL)(f) \right\|_{L^2(\R^{+} \times M)} \left\|(1-e^{-sL^*})sL^*e^{-sL^*} g\right\|_{L^2(\R^+\times M)},
\end{align*}
where $\R^+ \times M$ is equipped with the tensorial-product measure $\frac{ds}{s} \otimes d\mu(x)$. The second quantity is bounded by $\|g\|_{L^2}$ since the $L^2$-boundedness of the Littlewood-Paley function (see Remark \ref{rem:holo} applied to $L^*$). Since $T_{s^m}({\bf 1})$ is a Carleson measure, it is well-known that we get the following inequality, making appear a non-tangential maximal function
\be{eq:maxnt} \left\|\textrm{Main}(f)\right\|_{L^2} \lesssim \left\| \sup_{\genfrac{}{}{0pt}{}{(s,y)\in \R^+ \times M}{|x-y|\leq s}} |\phi(s^mL)f(y)| \right\|_{L^2}. \ee
This quantity can be estimated by $\|f\|_{L^2}$ since the maximal function is pointwisely bounded by ${\mathcal M}(f)$ thanks to the pointwise estimates on the heat kernel (\ref{eq:pointwise}) and Remark \ref{rem:holo}. \\
It remains also to study the ``Error term'' in (\ref{eq:error}). For this term, we use duality again as follows~:
$$ \left|\langle  \textrm{Error}(f),g \rangle \right| \leq \left| \int_{\epsilon}^R sLe^{-sL} (1-e^{-sL})(g)(x) \left[T_s \phi(sL) - T_s({\bf 1}) \phi(sL)\right](f)(x) \frac{ds}{s} dx \right|.$$
Cauchy-Schwarz inequality yields
\begin{align*}
\left|\langle  \textrm{Error}(f),g \rangle \right| & \\
& \lesssim \left( \int_{\epsilon}^R \int \left|sLe^{-sL}(1-e^{-sL})(g)(x)\right|^2 \frac{d\mu(x)dt}{t} \right)^{1/2} \\
& \ \ \left( \int \int_{\epsilon}^R \left|\left[T_s \phi(sL) - T_s({\bf 1}) \phi(sL)\right](f)(x) \right|^2 \frac{dsd\mu(x)}{s} \right)^{1/2}\\
& \lesssim \|g\|_{L^2} \left( \int \int_{\epsilon}^R \left|\left[T_s \phi(sL) - T_s({\bf 1}) \phi(sL)\right](f)(x) \right|^2 \frac{dsd\mu(x)}{s} \right)^{1/2},
\end{align*}
where we have used the $L^2$-boundedness of the square function (see Remark \ref{rem:holo}). Let fix $s\in(\epsilon,R)$ and consider
$$  I_s:= \int \left|\left[T_s \phi(sL) - T_s({\bf 1}) \phi(sL)\right](f)(x) \right|^2 d\mu(x) =  \int \left|T_s\left[\phi(sL)f - \phi(sL)f(x)\right](x) \right|^2 d\mu(x).$$
We choose $(Q_i)_i$ a bounded covering of the whole manifold $M$ with balls of radius $r:=s^{1/m}$, in order that
\begin{align*}
 I_s & \leq \sum_i \int_{Q_i} \left|T_s\left[\phi(sL)f - \phi(sL)f(x)\right](x) \right|^2 d\mu(x) \\
 & \leq \sum_i \int_{Q_i} \ \sup_{y\in Q_i} \left|T_s\left[\phi(sL)f - \phi(sL)f(y)\right](x) \right|^2 d\mu(x).
\end{align*} 
The supremum over $y\in Q_i$ can be estimated with the Sobolev inequality, Proposition \ref{prop:sobolev}. Let us fix the index $i$. For a large enough integer $M$, Proposition \ref{prop:sobolev} implies
\begin{align*}
 \sup_{y\in Q_i} \left|T_s\left[\phi(sL)f - \phi(sL)f(y)\right](x) \right|^2 & \\
& \hspace{-4cm} \lesssim \sum_{j\geq 0} 2^{-2\delta j} \frac{1}{\mu(2^j Q_i)} \int_{2^j Q_i} \left|(1+sL)^{M} \left[T_s\left[\phi(sL)f - \phi(sL)f(\cdot)\right](x)\right] (y)\right| ^2 d\mu(y), \end{align*}
with an other exponent $\delta$, which is still strictly bigger than $1$. \\
By expanding $(1+sL)^{M}$ and using $L({\bf 1})=0$, it comes 
\begin{align*}
(1+sL)^{M} \left[ T_s\left[\phi(sL)f - \phi(sL)f(\cdot)\right](x)\right] (y) =  & T_s\left[\phi(sL)f - \phi(sL)f(y)\right](x) \\
& - T_s({\bf 1}) (x) \left[(1+sL)^{M}-1\right]\phi(sL)f(y). 
\end{align*}
Consequently,
\begin{align*}
\left|T_s\left[\phi(sL)f - \phi(sL)f(y)\right](x) \right|^2  \lesssim & \sum_{j\geq 0} 2^{-2\delta j} \frac{1}{\mu(2^j Q_i)}\int_{2^j Q_i} \left|T_s\left[\phi(sL)f - \phi(sL)f(y)\right](x)\right|^2 d\mu(y)  \\
& \hspace{-3cm}  + \left|T_s({\bf 1}) (x)\right|^2 \sum_{j\geq 0} 2^{-2 \delta j} \frac{1}{\mu(2^j Q_i)}\int_{2^j Q_i}  \left|\left[(1+sL)^{M}-1\right]\phi(sL)f(y) \right| ^2 d\mu(y).
\end{align*}
Finally, we get
\begin{align*}
 I_s & \leq \sum_i \sum_{j\geq 0} 2^{-2 \delta j}\frac{1}{\mu(2^j Q_i)} \int_{Q_i} \int_{2^j Q_i} \left|T_s\left[\phi(sL)f - \phi(sL)f(y)\right](x)\right|^2 d\mu(y) d\mu(x)  \\
& \hspace{0.5cm} + \sum_i \sum_{j\geq 0} 2^{-2 \delta j}\frac{1}{\mu(2^j Q_i)} \left(\int_{Q_i}\left|T_s{\bf 1}\right|^2 d\mu\right) \left(\int_{2^j Q_i} \left|\left[(1+sL)^{M}-1\right]\phi(sL)f\right| ^2 d\mu \right). 
\end{align*}
Let us denote by $I_s^1$ and $I_s^2$ the two previous terms. \\
Concerning $I_s^1$, we use Fubini's Theorem and off-diagonal estimates (\ref{eq:prop2}) and deduce 
\begin{align*}
I_s^1 & \leq  \sum_{i} \sum_{j\geq 0} 2^{-2 \delta j} \frac{1}{\mu(2^j Q_i)} \int_{Q_i} \int_{2^j Q_i} \left|T_s\left[\phi(sL)f - \phi(sL)f(y)\right](x)\right|^2 d\mu(y) d\mu(x) \\
 & \lesssim  \sum_{j\geq 0} 2^{-2 \delta j} \sum_{\genfrac{}{}{0pt}{}{i,l_1}{Q_{l_1}\subset 2^j Q_i}} \int_{Q_i} \bigg[ \\
&  \sum_{l_2} \left(1+\frac{d(Q_{l_1},Q_{l_2})}{s^{1/m}}\right)^{-d-2N-\delta} \left( \int_{Q_{l_2}} \left|\phi(sL)f(x) - \phi(sL)f(y)\right|^2 \frac{d\mu(y)}{\mu(2^j Q_i)} \right)^{1/2}\bigg]^2 d\mu(x).
\end{align*}
To be precised, from (\ref{eq:prop2}) we have off-diagonal decays for $T_s(e^{-sL} \cdot)$ at the scale $s$. We can also obtain off-diagonal decays for $T_s(e^{-sL/2} \cdot)$ at the same scale $s$. So at the last line of the previous inequality, it should appear $e^{sL/2}\phi(sL)$ (and not just $\phi(sL)$). For an easy readibility, we prefer to keep the notation $\phi(sL)$~: indeed $\phi(sL)$ is equal to a polynomial term multiplied by $e^{-sL}$ ; so $\phi(sL)$ and $e^{sL/2}\phi(sL)$ satisfy the same off-diagonal estimates at the scale $s$. \\
Since $\delta>1$, with another constant $\delta'>1$ (indeed belonging in $(1,\delta)$), Cauchy-Schwartz inequality gives
\begin{align*}
\left[\sum_{l_2} \left(1+\frac{d(Q_{l_1},Q_{l_2})}{s^{1/m}}\right)^{-d-2N-\delta} \left(\int_{Q_{l_2}} \left|\phi(sL)f(x) - \phi(sL)f(y)\right|^2 \frac{d\mu(y)}{\mu(2^j Q_i)} \right)^{1/2}\right]^2  \\ 
\lesssim 
\left[\sum_{l_2} \left(1+\frac{d(Q_{l_1},Q_{l_2})}{s^{1/m}}\right)^{-d-2N-2\delta'}
\int_{Q_{l_2}} \left|\phi(sL)f(x) - \phi(sL)f(y)\right|^2 \frac{d\mu(y)}{\mu(2^j Q_i)} \right],
\end{align*}
since for each $l_1$
$$ \sum_{l_2} \left(1+\frac{d(Q_{l_1},Q_{l_2})}{s^{1/m}}\right)^{-d-2(\delta-\delta')} \lesssim 1.$$
Then Proposition \ref{prop:poincare} and Poincar\'e inequality $(P_2)$ yield
\begin{align*}
I_s^1 & \lesssim   \sum_{j\geq 0} 2^{-2 \delta' j} \sum_{\genfrac{}{}{0pt}{}{i,l_1,l_2}{Q_{l_1}\subset 2^j Q_i}} \left(1+\frac{d(Q_{l_1},Q_{l_2})}{s^{1/m}}\right)^{-d-2N-2\delta'} \left[\frac{\mu(Q_{i})}{\mu(2^jQ_i)}\int_{Q_{l_2}}  {\mathcal M}_{2-\epsilon}[\nabla\phi(sL)f]^2 d\mu \right. \\
& \hspace{1cm} \left. + \frac{\mu(Q_{l_2})}{\mu(2^jQ_i)} \int_{Q_i}  {\mathcal M}_{2-\epsilon}[\nabla\phi(sL)f]^2 d\mu\right]  \left( d(Q_{l_1},Q_{l_2})^2 + 2^{2j} s^{2/m}\right),
\end{align*}
which can be divided in two quantities (we have used that $d(Q_i,Q_{l_2})\leq d(Q_{l_1},Q_{l_2}) + 2^j s^{1/m}$). For the first one, we can compute the sum over $i$ since for every $l_1$
$$ \sum_{\genfrac{}{}{0pt}{}{i}{Q_{l_1}\subset 2^j Q_i}} \frac{\mu(Q_{i})}{\mu(2^jQ_i)} \lesssim  \sum_{\genfrac{}{}{0pt}{}{i}{Q_{i}\subset 2^j Q_{l_1}}} \frac{\mu(Q_{i})}{\mu(2^jQ_{l_1})} \lesssim 1,$$
where we used that $Q_{l_1}\subset 2^j Q_i$ implies $Q_{i}\subset 2^j Q_{l_1}$ and $2^{j+2}Q_{i} \simeq 2^{j+2} Q_{l_1}$. For the second one, we can estimate the sum over $l_2$ thanks to for every $l_1$
\begin{align*}
 \sum_{l_2} \left(1+\frac{d(Q_{l_1},Q_{l_2})}{s^{1/m}}\right)^{-d-2N-2(\delta'-1)} \mu(Q_{l_2}) &\lesssim \mu(Q_{l_1}) \sum_{l_2} \left(1+\frac{d(Q_{l_1},Q_{l_2})}{s^{1/m}}\right)^{-d-2(\delta'-1)} \\
 &\lesssim \mu(Q_{l_1}) \sum_{k} 2^{- k(d+2\delta'-2)} \sum_{\genfrac{}{}{0pt}{}{l_2}{d(Q_{l_2},Q_{l_1})\simeq 2^k s}}  1 \\
 &\lesssim \mu(Q_{l_1}) \sum_{k} 2^{-k(d+2\delta'-2)} 2^{kd} \lesssim \mu(Q_{l_1}),
\end{align*}
due to the doubling property, $\delta'>1$ and arguments similar to (\ref{eq:argument}). So we obtain
\begin{align*} 
I_s^1 & \lesssim  \sum_{j\geq 0} 2^{-2(\delta'-1) j} \sum_{l_1,l_2} \left(1+\frac{d(Q_{l_1},Q_{l_2})}{s^{1/m}}\right)^{-d-2N-2(\delta'-1)} \int_{Q_{l_2}}  {\mathcal M}_{2-\epsilon}[s^{1/m}\nabla\phi(sL)f]^2 d\mu \\
& \hspace{1cm}  +  
 \sum_{j\geq 0} 2^{-2(\delta'-1) j} \sum_{\genfrac{}{}{0pt}{}{i,l_1}{Q_{l_1}\subset 2^j Q_i}} \frac{\mu(Q_{l_1})}{\mu(2^jQ_i)} \int_{Q_i}  {\mathcal M}_{2-\epsilon}[s^{1/m}\nabla\phi(sL)f]^2 d\mu \\
 & \lesssim  \sum_{l_2} \int_{Q_{l_2}}  {\mathcal M}_{2-\epsilon}[s^{1/m} \nabla\phi(sL)f]^2 d\mu +
 \sum_{j\geq 0} 2^{-2(\delta'-1) j} \sum_{i} \int_{Q_i}  {\mathcal M}_{2-\epsilon}[s^{1/m}\nabla\phi(sL)f]^2 d\mu \\
 & \lesssim  \int {\mathcal M}_{2-\epsilon}[s^{1/m}\nabla\phi(sL)f]^2 d\mu.
\end{align*}
We have used the bounded-overlap property of $(Q_i)_i$ (balls of radius $s$). By the $L^2$-boundedness of the maximal operator ${\mathcal M}_{2-\epsilon}$ (see Theorem \ref{MIT}), we conclude to 
\begin{align}
I_s^1 \lesssim \int \left|s^{1/m}\nabla \phi(sL)f\right|^2 d\mu. \label{Is1}
\end{align}
It remains to study the second term $I_s^2$, which is equal to
\begin{align*}
 I_s^2 & = \sum_i \sum_{j\geq 0} 2^{-2\delta j}\frac{1}{\mu(2^j Q_i)} \left(\int_{Q_i}\left|T_s{\bf 1} (x)\right|^2 d\mu(x)\right) \left(\int_{2^j Q_i} \left|\left[(1+sL)^{M}-1\right]\phi(sL)f(y)\right| ^2 d\mu(y) \right).
\end{align*}
Since $T_s{\bf 1} (x) = (sL)^\K e^{-sL}(1-e^{-sL}) T({\bf 1})$ and $T({\bf 1})$ belongs to $BMO_L$, Proposition \ref{prop:Carleson} yields
$$ \int_{Q_i}\left|T_s{\bf 1} (x)\right|^2 dx \lesssim \mu(Q_i).$$
Hence, 
\begin{align}
 I_s^2 & \lesssim \sum_i \sum_{j\geq 0} 2^{-2\delta j}\frac{\mu(Q_i)}{\mu(2^j Q_i)} \int_{2^j Q_i} \left|\left[(1+sL)^{M}-1\right]\phi(sL)f(y)\right| ^2 d\mu(y)  \nonumber \\
 & \lesssim \int \left|\left[(1+sL)^{M}-1\right]\phi(sL)f\right| ^2 d\mu, \label{Is2}
\end{align}
where we used similar arguments as (\ref{eq:argument}) in order to prove that for each integer $j$, the collection $(2^j Q_i)_i$ is a $2^{jd}$-bounded covering~:
$$ \sum_{i} {\bf 1}_{2^j Q_i} \lesssim 2^{jd}.$$
Consequently, from (\ref{Is1}) and (\ref{Is2}) we obtain
$$I_s \lesssim \int \left|s^{1/m} \nabla \phi(sL)f\right|^2  + \left|\left[(1+sL)^{M}-1\right]\phi(sL)f\right|^2 d\mu,$$
which finally yields
\begin{align*}
 \lefteqn{\left|\langle  \textrm{Error}(f),g \rangle \right|} & & \\
& &  \lesssim \|g\|_{L^2} \left( \int_{\epsilon}^R\int \left|s^{1/m}\nabla \phi(sL)f(y)\right|^2 + \left|\left[(1+sL)^{M}-1\right]\phi(sL)f(y)\right|^2 \frac{d\mu(y) ds}{s} \right)^{1/2}.
\end{align*}
We also conclude to the desired estimate
$$ \|\textrm{Error}(f)\|_{L^2} \lesssim \|f\|_{L^2}$$
by invoking Assumption (\ref{ass:square}) for the first quantity and point 2 of Remark \ref{rem:holo} for the second one. Since the two terms $\textrm{Main}(f)$ and $\textrm{Error}(f)$ of (\ref{eq:error}) have been bounded in $L^2$, we have also proved that $U_{\epsilon,R}$ is uniformly (with respect to $\epsilon,R$) bounded on $L^2$. The second step is also finished and the proof of the first part of the theorem too. \findem

\gb
{\bf Proof of the second part of Theorem \ref{thm:principal}}. \\
We assume that $T$ admits a continuous extension from $L^2$ to $L^2$. By symmetry, it suffices us to check that $T({\bf 1})$ belongs to $BMO_L$. By definition of $BMO_L$ (see Subsection \ref{subsec:BMO}), we have to show that
$$ \sup_{\genfrac{}{}{0pt}{}{Q}{r_Q=t^{1/m}}} \frac{1}{\mu(Q)}\int_Q \left|T({\bf 1})-e^{-tL}T({\bf 1})\right| d\mu <\infty.$$
We recall that $T({\bf 1})$ is assumed to be well-defined. Indeed, we will prove a weaker but equivalent property
\be{BMO} \sup_{Q} \left(\frac{1}{\mu(Q)}\int_Q \left|T({\bf 1})-e^{-tL}T({\bf 1})\right|^2 d\mu\right)^{1/2} <\infty \ee
due to John-Nirenberg inequality (see Theorem 3.1 in \cite{DY} for a proof of such properties in this general framework concerning semigroups and \cite{BZ3} for an extension of such results in a more abstract setting). \\
Let us fix a ball $Q$ of radius $r:=t^{1/m}$ and by duality a function $f\in L^2(Q)$ such that $\|f\|_{L^2(Q)}\leq \mu(Q)^{-1/2}$. We have to bound
$$ \left|\langle (1-e^{-tL})  T({\bf 1}),f\rangle \right| = \left|\langle   T({\bf 1}), (1-e^{-tL^*})f\rangle \right|.$$
First we set $\chi_Q$ for the characteristic function of $4Q$ and by assumption, it comes
\begin{align*}
 \left|\langle   T(\chi_Q), (1-e^{-tL^*})f\rangle \right| &  \lesssim \|\chi_Q\|_{L^2} \|f\|_{L^2} \\
& \lesssim \mu(Q)^{1/2} \mu(Q)^{-1/2} \lesssim 1,
\end{align*}
where we used the doubling property and the $L^2$-boundedness of the semigroup. So we have to bound the remainder term, which we differentiate as follows
\begin{align*}
\left|\langle   T(1-\chi_Q), (1-e^{-tL^*})f\rangle \right| & \leq \int_0^t \left|\langle sLe^{-sL}T(1-\chi_Q),f\rangle \right| \frac{ds}{s}.
\end{align*}
For each $s$ and each integer $j\geq 1$, we consider $(Q^{s,j}_l)_l$ a bounded covering of $2^{j}Q\setminus 2^{j-1}Q$ by balls of radius $s^{1/m}$ and we associate $\chi^{s,j}_l$ a partition of unity, in order that
$$ \chi_Q + \sum_{j\geq 1} \sum_l \chi^{s,j}_l = {\bf 1}.$$
Similarly, let $(Q^{s}_l)_l$ be a bounded covering of $Q$ by balls of radius $s^{1/m}$ and we associate $\chi^{s}_l$ a partition of unity, in order that
$$ f = \sum_{l} f \chi^s_l .$$
From off-diagonal decay (\ref{assum:off}) and since $\K=1$, it yields
\begin{align*}
 \left|\langle sLe^{-sL}T(1-\chi_Q),f\rangle \right| & \lesssim \sum_{j\geq 0} \sum_{l_1,l_2} \left|\langle sLe^{-sL}T(\chi^{j,s}_{l_1}),f \chi^s_{l_2}\rangle \right| \\
 & \hspace{-2cm} \lesssim \sum_{j\geq 1} \sum_{l_1,l_2} \|\chi^{j,s}_{l_1}\|_{L^2} \|f\chi_{l_2}^s\|_{L^2} \left(1+\frac{2^jr}{s^{1/m}}\right)^{-d-\delta} \\
 & \hspace{-2cm} \lesssim \sum_{j\geq 1} \left(\sum_{l_1} \mu(Q^{j,s}_{l_1})\right)^{1/2} \left( \sum_{l_2} \|f\chi_{l_2}^s\|_{L^2}^2\right)^{1/2} \left(\frac{2^j r}{s^{1/m}}\right)^{d/2}\left(\frac{r}{s^{1/m}}\right)^{d/2} \left(\frac{s^{1/m}}{2^jr}\right)^{d+\delta} \\
 & \hspace{-2cm} \lesssim \sum_{j\geq 1} \mu(2^jQ)^{1/2} \mu(Q)^{-1/2} \left(\frac{r}{s^{1/m}}\right)^{-\delta} 2^{-j(d/2+\delta)} \\
 & \hspace{-2cm} \lesssim \sum_{j\geq 1} 2^{-j\delta} \left(\frac{s^{1/m}}{r}\right)^{\delta} \\
 & \hspace{-2cm} \lesssim\left(\frac{s^{1/m}}{r}\right)^{\delta}.
\end{align*}
We have used Cauchy-Schwartz inequality in $l_1$ and $l_2$ with as previously (see arguments employed for (\ref{eq:argument}))
$$ \sum_{l_1} 1 \lesssim \left(\frac{2^j r}{s^{1/m}}\right)^{d} \quad \textrm{and} \quad \sum_{l_2} 1 \lesssim \left(\frac{ r}{s^{1/m}}\right)^{d}.$$
By integrating for $s\in(0,t)$, we deduce 
\begin{align*}
\left|\langle   T(1-\chi_Q), (1-e^{-tL^*})f\rangle \right| & \leq \int_0^t \left(\frac{s^{1/m}}{r}\right)^{\delta} \frac{ds}{s}  \lesssim 1,
\end{align*} 
which concludes the proof of (\ref{BMO}). \findem

\section[Applications]{Applications to new paraproducts and to Calder\'on-Zygmund operators on a Riemannian manifold} \label{sec:para}
 
We consider an operator $L$ satisfying the assumptions of the previous sections. Moreover we will assume off-diagonal decays for the gradient of the semigroup, as follows: for all $s>0$, all balls $Q_1,Q_2$ of radius $r=s^{1/m}$, we have for every integer $k\geq 0$
\be{ass:gradient2} \left\| s^{1/m} \nabla (sL)^k e^{-sL} (f)\right\|_{L^2(Q_2)} \lesssim \left(1+\frac{d(Q_1,Q_2)}{r}\right)^{-d-2N-\delta} \|f\|_{L^2(Q_1)}. \ee
We just emphasize that this new assumption still holds in the examples described in Subsection \ref{subsec:ex}. Indeed usually, gradient of the semigroup satisfies Gaffney estimates.

\mb We want to study new kind of paraproducts, relying on the semigroup. For $b\in L^\infty \subset BMO_L$ (indeed since Remark 7.6 in \cite{BZ} and Proposition 6.7 in \cite{DY1}, we know that $L^\infty \subset BMO\subset BMO_{L}$ thanks to $L({\bf 1})=0$) and $f,g\in L^2$, let us consider the trilinear form
$$ \Lambda^1(b,f,g) :=  \int_M \int_{0}^\infty \left[\psi_t(L^*)g\right] \left[\phi_t(L)b \, \psi_t(L)f\right] \frac{dt}{t} d\mu,$$
where we write for convenience
$$ \psi_t(L):=(tL)^N e^{-tL}(1-e^{-tL}) \quad \textrm{and} \quad \phi_t(L):=e^{-tL},$$
with a large enough integer $N>d/m$.

\mb As a direct consequence of Cauchy-Schwartz inequality, pointwise bound on $\phi_t(L)$ and quadratic estimates (due to Remark \ref{rem:holo}), we know that $\Lambda^1$ is bounded on $L^\infty \times L^2 \times L^{2}$. Then we deduce the following result.

\begin{prop} \label{prop:para1} The trilinear form $\Lambda^1$ is bounded on $L^\infty \times L^p \times L^{p'}$ for every exponent $p\in (1,\infty)$.
\end{prop}

\dem Let us fix the function $b\in L^\infty$ and consider the linear operator $U$ such that
$$ \langle U(f),g \rangle := \Lambda^1(b,f,g).$$
It is given by
$$ U(f):= \int_0^\infty  \psi_t(L)\left[\phi_t(L)b\, \psi_t(L)f \right] \frac{dt}{t}.$$
We have just seen that $T$ is bounded on $L^2$. Then, we let the details to the reader and refer to Proposition \ref{thm:paraproduit} (where we prove a stronger result). It is quite easy to check that $U$ satisfies the assumptions
(\ref{assum:off}), (\ref{assum:off-dual}) and (\ref{assum:weaklybounded}). We also deduce the desired result by applying Corollary \ref{cor}. \findem 

\mb We are now looking to invert the role of the $L^\infty$-function $b$ and the $L^2$-function $f$~:

\begin{prop} \label{prop:para2} The trilinear form $\Lambda^2$ defined by
$$ \Lambda^2(b,f,g) := \int_M \int_{0}^\infty \left[\psi_t(L^*)g\right] \left[ \phi_t(L)f \,\psi_t(L)b\right] \frac{dt}{t} d\mu,$$
is bounded on $L^\infty \times L^p \times L^{p'}$ for every exponent $p\in (1,\infty)$.
\end{prop}

\dem Using the Carleson measure property (Proposition \ref{prop:Carleson} and since $L^\infty\subset BMO_L$) together with Cauchy-Schwartz inequality, we obtain the desired result for $p=2$. By the same reasoning as used for Proposition \ref{prop:para1}, we conclude this proof.
\findem

\mb Using tri-linear interpolation and symmetry, we deduce the following result.

\begin{prop} \label{prop:para3} The trilinear form $\Lambda$ defined by
$$ \Lambda(h,f,g) := \int_{0}^\infty \int_M \psi_t(L)g \, \phi_t(L)f \, \psi_t(L)h \frac{dt}{t} d\mu,$$
is bounded on $L^p \times L^q \times L^{r}$ for every exponents $p,q,r\in (1,\infty]$ satisfying
$$ \frac{1}{p}+\frac{1}{q}+\frac{1}{r} = 1.$$
\end{prop}

\mb These results concerning paraproducts with two functions $\psi_t$ are also easily obtained, thanks to duality and Cauchy-Schwartz inequality in the variable $t$. Let us note that we do not need $N>d/m$. We are now interesting to paraproducts, involving only one function $\psi_t$.

\begin{rem} \label{rem:r} Let us first examine to this situation in the ``classical case''. Let us consider the Euclidean space $M=\R^d$ and denote by $\Psi$ a smooth function on $\R^d$ whose its spectrum is contained in a corona around $0$ and $\Phi$ another smooth function with a bounded spectrum. Then with the usual notations $\Psi_t:=t^{-d} \Psi(./t)$ and similarly for $\Phi$, we are interested in the following paraproducts
$$ f \rightarrow \int_0^\infty \Psi_t\left[ \Phi_t(f) \Phi_t(b) \right] \frac{dt}{t},$$
with $b\in L^\infty$. By duality, it gives rise to the following trilinear function~:
$$ (b,f,g) \rightarrow \int_0^\infty\int_{\R^d} \Psi_t(g) \Phi_t(f) \Phi_t(b) \frac{dtdx}{t}.$$
Since we know that the spectrum of a product is contained in the Minkowski sum of the spectrums, it follows that we can find some ``good'' smooth functions $\tilde{\Psi}$ and $\tilde{\Phi}$ (satisfying the same spectral property than $\Psi$ and $\Phi$) such that
\begin{align*}
  \int_{\R^d} \int_0^\infty \Psi_t(g) \Phi_t(b) \Phi_t(b) \frac{dtdx}{t} = & \int_{\R^d} \int_0^\infty \Psi_t(g) \tilde{\Phi}_t(f) \tilde{\Psi}_t(b) \frac{dtdx}{t} \\
& + \int_{\R^d} \int_0^\infty \Psi_t(g) \tilde{\Psi}_t(f) \tilde{\Phi}_t(b) \frac{dtdx}{t}.
\end{align*}
So it comes that such paraproducts can be reduced to the sum of two paraproducts involving ``two functions $\Psi$''. Hence, they are bounded in Lebesgue spaces. \\ 
This reduction is due to the frequency analysis of the product. It is not clear how we can apply a similar reasoning in the framework of semigroup. This is the goal of the two following subsections, via our T(1) theorem.
\end{rem}

\subsection{Boundedness of new paraproducts in Lebesgue spaces with $r'>1$}

\begin{thm} \label{thm:paraproduit} The trilinear form $\Lambda$ defined by
$$ \Lambda(h,f,g) := \int_M \int_{0}^\infty \left[\psi_t(L^*)g\right] \left[ \phi_t(L)f \, \phi_t(L)h\right] \frac{dt}{t} d\mu,$$
is bounded on $L^p \times L^q \times L^{r}$ for every exponents $p,q,r\in (1,\infty)$ satisfying
$$ \frac{1}{p}+\frac{1}{q}+\frac{1}{r} = 1.$$
Equivalently, the paraproduct
$$ (h,f) \rightarrow \int_{0}^\infty \psi_t(L)\left[\phi_t(L)f \, \phi_t(L)h\right] \frac{dt}{t}$$
is bounded from $L^p \times L^q$ to $L^{r'}$. \\
Moreover $p$ or $q$ may be infinite.
\end{thm}

\dem By trilinear interpolation, it suffices us to prove boundedness for the limiting case: when one of the exponents is infinite. By symmetry between $f$ and $h$, we have also to deal with only one case: when $p=\infty$ (step 1) and then conclude by interpolation (step 2).

\mb
{\bf Step 1:} Estimate for $p=\infty$. \\
Let us fix $h\in L^\infty$ and consider the operator $U$ defined by
$$ U(f):= \int_0^\infty \psi_t(L)\left[\phi_t(L)f \, \phi_t(L)h\right] \frac{dt}{t}$$
in order that
$$ \Lambda(h,f,g):=\langle U(f),g\rangle.$$
We will prove that $U$ satisfies Assumptions of Theorem \ref{thm:principal}, then the desired result will follow from Corollary \ref{cor}. First, $L({\bf 1})=0$ yields
$$ U({\bf 1}) = \int_0^\infty \psi_t(L)\phi_t(L) h \frac{dt}{t} = ch\in L^\infty \subset BMO_L$$
where $c:=\int_0^\infty \psi_t(x) \phi_t(x) \frac{dt}{t}$ is a numerical constant independent on $x$. Moreover, it is obvious that 
$$ U^*({\bf 1}) = 0.$$
So it remains to check Assumptions (\ref{assum:off}), (\ref{assum:off-dual}), (\ref{assum:weaklybounded}), which we recall here. \\
Consider a large enough integer $\K$. For balls $Q_1,Q_2$ of radius $r=s^{1/m}$, we have
\begin{itemize}
\item if $d(Q_1,Q_2)\geq 2 r$, then
\be{assum:offbis} \left\|(sL)^\K e^{-sL} U(f)\right\|_{L^2(Q_2)} \lesssim \left(1+\frac{d(Q_1,Q_2)}{r}\right)^{-d-2N-\delta} \|f\|_{L^2(Q_1)}. \ee
and the dual estimates
\be{assum:off-dualbis} \left\|(sL^*)^\K e^{-sL^*} U^*(f)\right\|_{L^2(Q_2)} \lesssim \left(1+\frac{d(Q_1,Q_2)}{r}\right)^{-d-2N-\delta} \|f\|_{L^2(Q_1)}. \ee
\item if $d(Q_1,Q_2)\leq 2 r$, then we have
\be{assum:weaklyboundedbis} \left\|(sL)^\K e^{-sL} U(e^{-sL}f)\right\|_{L^2(Q_2)} + \left\|(sL^*)^{\K} e^{-sL^*} U(e^{-sL^*}f)\right\|_{L^2(Q_2)} 
\lesssim \|f\|_{L^2(Q_1)}. \ee
\end{itemize}

\mb
{\bf Step 1-1:} Assumption (\ref{assum:offbis}). \\
The operator $U$ is given by
$$ U(f):= \int_0^\infty \psi_t(L)\left[\phi_t(L)f \, \phi_t(L)h\right] \frac{dt}{t}.$$
We divide the integral on $t$ as follows
$$
\left(\int_{Q_2} \left|(sL)^\K e^{-sL} U(f)\right|^2 d\mu\right)^{1/2} \lesssim I+II$$
with 
$$ I:= \int_0^s \left\| s^\K t^N L^{\K+N} e^{-(s+t)L}(1-e^{-tL}) \left[\phi_t(L)f \, \phi_t(L)h\right] \right\|_{L^2(Q_2)} \frac{dt}{t} $$
and
$$ II:= \int_s^\infty \left\| s^\K t^N L^{\K+N} e^{-(s+t)L}(1-e^{-tL}) \left[\phi_t(L)f \, \phi_t(L)h\right] \right\|_{L^2(Q_2)} \frac{dt}{t}.$$
Let us treat the first term $I$ (the reasoning is similar to the one used for Corollary \ref{cor}). \\
Thanks to the $L^2$-off diagonal decay of the semigroup (and its derivative), we have since $t+s\in[s,2s]$
\begin{align*}
I & \lesssim  \int_0^s \left\| s^\K t^N L^{\K+N} e^{-(s+t)L}(1-e^{-tL}) \left[\phi_t(L)f \, \phi_t(L)h\right] \right\|_{L^2(Q_2)} \frac{dt}{t} \\
& \lesssim \int_0^s \frac{s^\K t^N}{(s+t)^{\K+N}}\left\| (s+t)^{\K+N}L^{\K+N} e^{-(s+t)L}(1-e^{-tL}) \left[\phi_t(L)f \, \phi_t(L)h\right] \right\|_{L^2(Q_2)} \frac{dt}{t}.
\end{align*}
Let now fix $t\in (0,s]$. Then $s+t\in [s,2s]$ so we know that $(s+t)^{\K+N} L^{\K+N} e^{-(s+t)L}$ satisfies $L^2-L^2$ off-diagonal decay at the scale $s$. Hence by considering $(R_k)_k$ a covering of the whole space $M$ by balls of radius $r$, we get
\begin{align*}
\left\| (s+t)^{\K+N} L^{\K+N} e^{-(s+t)L}(1-e^{-tL}) \left[\phi_t(L)f \, \phi_t(L)h\right] \right\|_{L^2(Q_2)} &  \\
& \hspace{-7cm} \lesssim \sum_{k} \left(1+\frac{d(Q_2,R_k)}{r}\right)^{-d-2N-\delta} \| (1-e^{-tL}) \left[\phi_t(L)f \,  \phi_t(L)h\right]\|_{L^2(R_k)}.
\end{align*}
Since $(1-e^{-tL})$ and $\phi_t(L)$ satisfy $L^2-L^2$ off-diagonal decays at the scale $t$ and $\phi_t(L)h$ is pointwisely bounded, it comes that $f\to (1-e^{-tL})\left[\phi_t(L)f \, \phi_t(L)h\right]$ satisfies similar $L^2-L^2$ off-diagonal decays.
So let $(R_k^j)_j$ a bounded covering of $R_k$ by balls of radius $t^{1/m}$ (and similarly for $Q_1$), it yields
\begin{align*}
\| (1-e^{-tL}) \left[\phi_t(L)f \, \phi_t(L)h\right]\|_{L^2(R_k)} & \lesssim \left(\sum_{j_1} \| (1-e^{-tL}) \left[\phi_t(L)f \, \phi_t(L)h\right]\|_{L^2(R_k^{j_1})}^2\right)^{1/2} \\
& \hspace{-3cm} \lesssim \left(\sum_{j_1} \left( \sum_{j_2} \left(1+\frac{d(R_k^{j_1},Q_1^{j_2})}{t^{1/m}}\right)^{-d-2N-\delta}  \|f\|_{L^2(Q_1^{j_2})}\right)^2\right)^{1/2} \\
& \hspace{-3cm} \lesssim \left(1+\frac{d(R_k,Q_1)}{t^{1/m}}\right)^{-d-2N-\delta} \left(\sum_{j_1} \left( \sum_{j_2}   \|f\|_{L^2(Q_1^{j_2})}\right)^2\right)^{1/2} \\
& \hspace{-3cm} \lesssim \left(1+\frac{d(R_k,Q_1)}{t^{1/m}}\right)^{-d-2N-\delta} \|f\|_{L^2(Q_1)}\left(\sum_{j_1,j_2}  1 \right)^{1/2} \\
& \hspace{-3cm} \lesssim \left(1+\frac{d(R_k,Q_1)}{t^{1/m}}\right)^{-d-2N-\delta} \|f\|_{L^2(Q_1)}\left(\frac{s}{t} \right)^{d/m},
\end{align*}
where we refer the reader to (\ref{eq:argument}) for the estimate of the sum over $j_1,j_2$. Finally, we get
\begin{align*}
\left\| (s+t)^{\K+N} L^{\K+N} e^{-(s+t)L}(1-e^{-tL}) \left[\phi_t(L)f \, \phi_t(L)h\right] \right\|_{L^2(Q_2)} &  \\
& \hspace{-7cm} \lesssim \sum_{k} \left(1+\frac{d(Q_2,R_k)}{s^{1/m}}\right)^{-d-2N-\delta} \left(1+\frac{d(R_k,Q_1)}{t^{1/m}}\right)^{-d-2N-\delta} \|f\|_{L^2(Q_1)}\left(\frac{s}{t} \right)^{d/m} \\
& \hspace{-7cm} \lesssim  \left(1+\frac{d(Q_2,Q_1)}{r}\right)^{-d-2N-\delta} \left(\frac{s}{t} \right)^{d/m} \|f\|_{L^2(Q_1)}.
\end{align*}
This permits to deduce (\ref{assum:offbis}) for the first term $I$ since $N>d/m$ and
$$ \int_0^s \frac{s^\K t^N}{(s+t)^{\K+N}} \left(\frac{s}{t} \right)^{d/m} \frac{dt}{t}\lesssim  1.$$
Concerning the second term $II$, we produce a similar reasoning and we deduce that for $t\geq s$
$$ \left\| (s+t)^{\K+N} L^{\K+N} e^{-(s+t)L}(1-e^{-tL}) \left[\phi_t(L)f \, \phi_t(L)h\right] \right\|_{L^2(Q_2)} \lesssim \left(1+\frac{d(Q_2,Q_1)}{t^{1/m}}\right)^{-d-2N-\delta} \|f\|_{L^2(Q_1)}.$$
Indeed, each appearing operator admits off-diagonal decay at the scale $t^{1/m}$ since $s+t\simeq t$. We also conclude to (\ref{assum:offbis}) for the first term $II$ since for $\K\geq d+2N+\delta$
$$ \int_s^\infty \frac{s^\K t^N}{(s+t)^{\K+N}} \left(1+\frac{d(Q_2,Q_1)}{t^{1/m}}\right)^{-d-2N-\delta} \frac{dt}{t}\lesssim \left(1+\frac{d(Q_2,Q_1)}{r}\right)^{-d-2N-\delta}.$$
We have also finished to check Assumption (\ref{assum:offbis}).

\mb
{\bf Step 1-2:} Assumption (\ref{assum:off-dualbis}). \\
By duality, (\ref{assum:off-dualbis}) is equivalent to
\be{assum:offbis2} \left\| U( (sL)^\K e^{-sL} f)\right\|_{L^2(Q_2)} \lesssim \left(1+\frac{d(Q_1,Q_2)}{r}\right)^{-d-2N-\delta} \|f\|_{L^2(Q_1)}, \ee
which we are going to prove. \\
We will just give the sketch of the proof since technical details are by now routine. As previously, the quantity to estimate can be divided in two quantities $I^*$ and $II^*$ with
$$ I^*:= \left\| \int_0^s (t L)^{N} e^{-tL}(1-e^{-tL}) \left[(sL)^\K e^{-sL}\phi_t(L)f \, \phi_t(L)h\right] \frac{dt}{t} \right\|_{L^2(Q_2)}  $$
and
$$ II^*:= \int_s^\infty \left\|  (t L)^{N} e^{-tL}(1-e^{-tL}) \left[(sL)^\K e^{-sL}\phi_t(L)f \, \phi_t(L)h\right] \right\|_{L^2(Q_2)} \frac{dt}{t}.$$
The second one can be exactly estimated as $II$ in the previous point and so we only deal with the first one. First, it is easy to check that we can replace the above quantity
$$  (tL)^{N} e^{-tL}(1-e^{-tL}) \left[(sL)^\K e^{-sL}\phi_t(L)f \, \phi_t(L)h\right]$$ by
$$  (tL)^{N} e^{-tL}(1-e^{-tL}) \left[\tilde{\psi}_s(L)f \, \phi_t(L)h\right],$$
where $\tilde{\psi}_s(L) = (sL)^\K e^{-sL}$ since $ (sL)^\K e^{-sL}\phi_t(L) = (sL)^\K e^{-(s+t)L}$ and $t\in(0,s]$. Indeed by computing the difference, it appears  
$$ e^{-sL}-e^{-(s+t)L} \simeq t Le^{-sL} = \frac{t}{s} (sL) e^{-sL},$$ 
involving an extra factor $\frac{t}{s}$ which permits to easily bound the difference as desired. So let us just consider
\be{eq:1}  \left\| \int_0^s  (t L)^{N} e^{-tL}(1-e^{-tL}) \left[\tilde{\psi}_s(L)f \, \phi_t(L)h\right] \frac{dt}{t} \right\|_{L^2(Q_2)}.\ee 
Since $\tilde{\psi}_s(L)f$ is essentially constant at the scale $t<<s$, we can compare the previous quantity to the following one
\be{eq:2}  \left\|  \tilde{\psi}_s(L)f  \int_0^s  (t L)^{N} e^{-tL}(1-e^{-tL}) \left[\phi_t(L)h\right] \frac{dt}{t}  \right\|_{L^2(Q_2)},\ee 
which satisfies the desired estimate since $\int_0^s  (t L)^{N} e^{-tL}(1-e^{-tL}) \left[\phi_t(L)h\right] \frac{dt}{t}$ is a uniformly bounded function for $h\in L^\infty$. It also remains us to study the difference between (\ref{eq:1}) and (\ref{eq:2}). We let to the reader the details, but the analysis of the difference is based on exactly the same arguments as used for the study of $I_s$, in the proof of Theorem \ref{thm:principal}. The difference makes appear the gradient $\nabla \tilde{\psi}_s(L)f$ at the scale $t$, so it comes an extra decay like
$$ t^{1/m}\nabla \tilde{\psi}_s(L)f = \left(\frac{t}{s}\right)^{1/m} s^{1/m}\nabla \tilde{\psi}_s(L)f.$$
Since Assumption (\ref{ass:gradient2}), we obtain the desired off-diagonal decays and the extra factor $\left(\frac{t}{s}\right)^{1/m}$ permits one more time to make the integral on $t$ convergent.

\mb
{\bf Step 1-3:} Assumption (\ref{assum:weaklyboundedbis}). \\
We let to the reader to check that the two previous points (Steps 1-1 and 1-2) still holds when $d(Q_1,Q_2)\leq r$ and permit to prove   (\ref{assum:weaklyboundedbis}). 

\mb This finishes the proof of the step 1 and by Corollary \ref{cor}, it yields the desired estimates for $p=\infty$. 

\mb
{\bf Step 2:} End of the proof. \\
By symmetry between $f$ and $h$, we know that the trilinear form $\Lambda$ is bounded on $L^\infty \times L^q \times L^{q'}$ (step 1) and on $L^p \times L^\infty \times L^{p'}$ (by symmetry). So for $r\in (1,\infty)$ fixed, we know that $\Lambda$ is bounded on $L^\infty \times L^{r'} \times L^r$ and on $L^{r'} \times L^\infty \times L^r$, which by bilinear interpolation gives a boundedness on $L^p \times L^q \times L^r$.
\findem

\mb By duality, we have the following results~:

\begin{thm} \label{thm:paraproduit2} The trilinear form $\Lambda$ defined by
$$ \Lambda(h,f,g) := \int_{0}^\infty \int_M \left[\phi_t(L^*)g\right] \left[ \psi_t(L)f \, \phi_t(L)h\right] \frac{dt}{t} d\mu,$$
is bounded on $L^p \times L^q \times L^{r}$ for every exponents $p,q,r\in (1,\infty)$ satisfying
$$ \frac{1}{p}+\frac{1}{q}+\frac{1}{r} = 1.$$
Equivalently, the paraproduct
$$ (h,f) \rightarrow \int_{0}^\infty \phi_t(L)\left[\psi_t(L)f \, \phi_t(L)h\right] \frac{dt}{t}$$
is bounded from $L^p \times L^q$ to $L^{r'}$. \\
Moreover $p$ or $q$ may be infinite.
\end{thm}

\begin{rem} The tri-linear forms of Theorems \ref{thm:paraproduit} and \ref{thm:paraproduit2} naturally appear, in the study of the product. For example, let $\phi$ be the function
$$ \phi(x):=-\int_x^\infty y e^{-y}(1-e^{-y}) dy.$$
We let to the reader to check that all the previous results still hold with the new operator
$$ \phi_t(L):=\phi(tL)$$
(instead of $\phi_t(L)=e^{-tL}$). Then we get a ``spectral'' decomposition of the identity as follows~: up to some numerical constant $c$, we have
$$ f  = c \int_0^\infty \phi'(tL) f \frac{dt}{t}$$
according to Remark \ref{rem:holo}. So for two smooth functions, we have
$$ fg := c^3 \int_{s,u,v>0} \phi'(sL) \left[\phi'(uL)f \, \phi'(vL)g \right] \frac{dsdudv}{suv}.$$
Since $\phi'(x)=\psi(x):=xe^{-x}(1-e^{-x})$, it comes that (by integrating according to $t:=\min\{s,u,v\}$)
\begin{align*}
 fg := & c^3 \int_{0}^\infty \psi(tL) \left[\phi(tL)f \, \phi(tL)g \right] \frac{dt}{t} + c^3 \int_{0}^\infty \phi(tL) \left[\psi(tL)f \, \phi(tL)g \right] \frac{dt}{t} \\
& c^3\int_{0}^\infty \phi(tL) \left[\phi(tL)f \, \psi(tL)g \right] \frac{dt}{t}.
\end{align*}
Consequently, the pointwise product $fg$ can be decomposed with tree paraproducts (involving only one function $\psi$), studied by Theorems \ref{thm:paraproduit} and \ref{thm:paraproduit2}.
\end{rem}

\subsection{Boundedness in weighted Lebesgue spaces and extrapolation to exponents $r'\leq 1$}

We are now interesting to extend previous results to weighted Lebesgue spaces and with exponents $r'\leq 1$. We move the reader to Definition \ref{def:poids} for the usual class of weights.

\begin{prop} \label{prop:weight} Let $\psi$ and $\phi$ be defined as in Theorems \ref{thm:paraproduit} and \ref{thm:paraproduit2}.
Let $p,q,r\in (1,\infty)$ be exponents satisfying
$$ \frac{1}{p}+\frac{1}{q}+\frac{1}{r} = 1$$
and consider a weight $\omega$ belonging to ${\mathbb A}_{p}\cap {\mathbb A}_{q}$.
The paraproducts
$$ (h,f) \rightarrow \int_{0}^\infty \psi_t(L)\left[\phi_t(L)f \, \phi_t(L)h\right] \frac{dt}{t}$$
and
$$ (h,f) \rightarrow \int_{0}^\infty \phi_t(L)\left[\psi_t(L)f \, \phi_t(L)h\right] \frac{dt}{t} $$
are bounded from $L^p(\omega) \times L^q(\omega)$ to $L^{r'}(\omega)$.
\end{prop}

\dem Theorems \ref{thm:paraproduit} and \ref{thm:paraproduit2} corresponds to the desired result with the constant weight $\omega={\bf 1}$. \\
It is well-known that weighted estimates are closely related to estimates of some maximal sharp functions.
Let us denote the following maximal sharp function associated to an exponent $s\in [1,\infty)$ (introduced by J.M. Martell in \cite{Martell} and extended in \cite{BZ})
$$ M_s^\sharp (h)(x) := \left(\sup_{t>0} \frac{1}{B(x,t^{1/m})} \int_{B(x,t^{1/m})} \left| \psi(tL)(h) \right|^s d\mu \right)^{1/s}.$$
Let us explain how can we obtain the desired result only for the first paraproduct (the reasoning for the second one beeing similar)
$$ T(h,f):=\int_{0}^\infty \psi_t(L)\left[\phi_t(L)f \, \phi_t(L)h\right] \frac{dt}{t}.$$
The $L^2-L^2$ off-diagonal decays (\ref{assum:offbis}),  (\ref{assum:off-dualbis}) and (\ref{assum:weaklyboundedbis}) yield that
$$ M_2^\sharp (T(h,f))(x) \lesssim {\mathcal M}_2(f)(x) {\mathcal M}_2 (h)(x).$$
We let to the reader to check this point, but it is a direct consequence of the off-diagonal decays and the pointwise bound of the heat kernel (see Theorem 6.1 in \cite{BZ} for a detailed proof of such inequalities). In addition since we assume $L^1-L^\infty$ off-diagonal decays of the semigroup and its derivatives (pointwise estimates of the heat kernel), it is easy to see that we can obtain off-diagonal decays (\ref{assum:offbis}),  (\ref{assum:off-dualbis}) and (\ref{assum:weaklyboundedbis}) for all exponents $s\in (1,\infty)$ and not only for $s=2$. So we can obtain an estimate like
\be{eq:prod} M_s^\sharp (T(h,f))(x) \lesssim {\mathcal M}_s(f)(x) {\mathcal M}_s (h)(x) \ee
for every exponent $s>1$.\\
In order to compare the Lebesgue norm of $T(h,f)$ and the one of $M_s^\sharp (T(h,f))$, we need to use a Fefferman-Stein inequality. We refer the reader to \cite{Martell} (Theorem 4.2), \cite{BZ} (Corollary 5.8 and Theorem 6.4) and to \cite{B} (Lemma 2, Remark 3 and Lemma 3 for the weighted version) for such inequalities. Since pointwise estimates on the semigroup and Step 3 in the proof of Corollary \ref{cor}, we know that for every $s\in (1,\infty)$ and weight $\nu$
 $$ \left\| {\mathcal M}_s [T(h,f)] \right\|_{L^r(\nu)} \lesssim  \left\| M_s^\sharp [T(h,f)] \right\|_{L^r(\nu)}.$$
Such inequalities are based on ``good-$\lambda$ inequalities'' relatively to the two maximal operators, obtained in a very general framework by P. Auscher and J.M. Martell in \cite{AM} (Theorem 3.1). \\
Consequently, for our weight $\omega$, H\"older inequality and (\ref{eq:prod}) give
\begin{align*}
  \left\| T(h,f) \right\|_{L^r(\omega)}& \leq \left\| {\mathcal M}_s [T(h,f)] \right\|_{L^r(\omega)} \\
 & \lesssim \left\| M_s^\sharp [T(h,f)] \right\|_{L^r(\omega)} \\
& \lesssim \left\| {\mathcal M}_s(f) {\mathcal M}_s (h) \right\|_{L^r(\omega)} \\
& \lesssim \left\| {\mathcal M}_s(f)\right\|_{L^q(\omega)} \left\| {\mathcal M}_s (h) \right\|_{L^p(\omega)}.
\end{align*}
Then we chose $s\in(1,\min\{p,q\})$ in order that ${\mathcal M}_s$ is bounded in $L^p(\omega)$ and in $L^q(\omega)$ (due to $\omega \in {\mathbb A}_{p} \cap {\mathbb A}_{q}$) and so we conclude that $T$ is bounded from $L^p(\omega) \times L^q(\omega)$ into $L^{r'}(\omega)$.
\findem

\begin{rem} Using recent works of L. Grafakos, L. Liu, A. Lerner, S. Ombrosi, C. P\'erez, R. H. Torres and R.
 Trujillo-Gonz\'alez \cite{LOPTT, GLPT} ; it seems possible to get similar results with different weights for $h$ and $f$. This requires the notion of bilinear ${\mathbb A}_{\overrightarrow{P}}$ condition and the use of a ``bilinear strong maximal function''. We do not detail these possible improvements here.
\end{rem}

\mb Then we use theory of extrapolation to obtain new boundedness for our paraproducts (see Theorem 2 of \cite{GM}):

\begin{thm} \label{thm:paraproduit3} Let $\psi$ and $\phi$ be defined as in Theorems \ref{thm:paraproduit} and \ref{thm:paraproduit2}.
Let $p,q\in (1,\infty)$ and $r'\in(1/2,\infty)$ be exponents satisfying
$$ \frac{1}{p}+\frac{1}{q} = \frac{1}{r'}$$
and consider a weight $\omega$ belonging to ${\mathbb A}_{p}\cap {\mathbb A}_{q}$.
The paraproducts
$$ (h,f) \rightarrow \int_{0}^\infty \psi_t(L)\left[\phi_t(L)f \, \phi_t(L)h\right] \frac{dt}{t}$$
and
$$ (h,f) \rightarrow \int_{0}^\infty \phi_t(L)\left[\psi_t(L)f \, \phi_t(L)h\right] \frac{dt}{t} $$
are bounded from $L^p(\omega) \times L^q(\omega)$ to $L^{r'}(\omega)$.
\end{thm}

\begin{rem} The improvement in this new result is that the exponent $r'$ could be smaller than one.
\end{rem}

\subsection{A ``classical'' $T(1)$ Theorem on a Riemannian manifold} \label{subsec:manifold}

We devote this subsection to the proof of a T(1) theorem for Calder\'on-Zygmund operators on a general doubling Riemannian manifold $(M,d,\mu)$ of infinite measure.

\begin{df} \label{def:CZ} A function $K$ defined on $M \times M \setminus \{(x,x),x\in M \}$ is called a ``standard kernel of order $1+\epsilon$'' for some $\epsilon \in (0,1]$ if for all $x\neq y$
$$ \left|K(x,y)\right| \lesssim \frac{1}{d(x,y)}, $$ 
for $x'\in M$ satisfying $|x-x'|\leq \frac{1}{2} \max\{ |x-y|,|x'-y|\}$
$$ \left|\nabla K(x,y)- \nabla K(x',y)\right| \lesssim \frac{d(x,x')^\epsilon}{(d(x,y)+d(x',y))^{d+3N+1+\epsilon}} $$ 
and for $y'\in M$ satisfying $|y-y'|\leq \frac{1}{2} \max\{ |x-y|,|x-y'|\}$
$$ \left|\nabla K(x,y)- \nabla K(x,y')\right| \lesssim \frac{d(y,y')^\epsilon}{(d(x,y)+d(x,y'))^{d+3N+1+\epsilon}}.$$ 
A linear operator $T$, continuously acting from ${\mathcal M}$ to ${\mathcal S}'(M)$ and satisfying the integral representation
$$ \forall f\in C^\infty_0(M),\ \forall x \notin \textrm{supp}(f) \qquad T(f)(x) = \int_{M} K(x,y) f(y) d\mu(y), $$
is said to be associated to the kernel $K$.
\end{df}

\begin{thm} Let us assume that the doubling manifold $M$ satisfies Poincar\'e $(P_2)$ and Assumption (\ref{assm}), on its heat kernel. For $T$ a linear operator associated to a standard kernel of order $\epsilon$ (such that $T({\bf 1})$ and $T^*({\bf 1})$ are well-defined in ${\mathbb M}$), the two following properties are equivalent~:
\begin{itemize}
 \item $T$ is bounded on $L^2$
 \item $T({\bf 1})$ and $T^*({\bf 1})$ belong to $BMO_{-\Delta}$ and $T$ satisfies to the weak-boundedness property (\ref{assum:weaklybounded}).
\end{itemize}
Since $BMO \subset BMO_{-\Delta}$, if $T({\bf 1})$ and $T^*({\bf 1})$ belong to $BMO$ then they belong to $BMO_{-\Delta}$ too.
\end{thm}

\dem We look for applying our new T(1) theorem as follows. Let us consider $L=-\Delta$ given by the Laplacian. Then we have seen if Subsubsection \ref{subsub:manifold} that under (\ref{assm}) all our required assumptions are satisfied by the heat semigroup $(e^{-tL})_{t>0}$ with $m=2$. Moreover, we know that we have pointwise gaussian bound of the heat kernel $p_t$~:
$$ \left|p_t(x,y)\right| \lesssim e^{-\gamma d(x,y)^2/t}$$
for some constant $\gamma >0$. \\
Since $BMO$ is included in $BMO_L$ (see Proposition 6.7 in \cite{DY1} and Remark 7.6 of \cite{BZ}), it remains us to check that our operator $T$ satisfies to (\ref{assum:off}), (\ref{assum:off-dual}) with $\K=1$. By duality and symmetry, we only deal with (\ref{assum:off-dual})~: for every $s>0$, every ball $Q_1,Q_2$ of radius $r:=s^{1/2}$ (with $d(Q_1,Q_2)\geq 2r$) and function $f\in L^2(Q_1)$ 
\be{assum:offfin} \left\| (-s\Delta) e^{s\Delta}T(f)\right\|_{L^2(Q_2)} \lesssim \left(1+\frac{d(Q_1,Q_2)}{r}\right)^{-d-2N-\delta} \|f\|_{L^2(Q_1)}. \ee
So let us consider the balls $Q_1$ and $Q_2$ and write for $x_0\in Q_2$
\begin{align*}
(-s\Delta) e^{s\Delta} T(f) (x_0) & = -\int s \Delta_x  p_s(x_0,y) T(f)(y) d\mu(y) \\
& = -\int s \Delta_y  p_s(x_0,y) T(f)(y) d\mu(y) \\
& = \int s \nabla_y  p_s(x_0,y)  \nabla T(f)(y) d\mu(y) \\
& = \int s \nabla_y  p_s(x_0,y) \left[\nabla T(f)(y) - \nabla T(f)(x_0)\right] d\mu(y),
\end{align*}
where we used the self-adjoint properties of the Laplacian, an integration by parts and at the last line the fact that
$$ \int \nabla_y  p_s(x_0,y) d\mu(y) = 0.$$
So it comes
\begin{align*}
\left|(-s\Delta) e^{s\Delta} T(f) (x_0)\right| & \leq \int\int s \left|\nabla_y  p_s(x_0,y)\right| \left| \nabla_y K(y,z) - \nabla_y K(x_0,z)\right| |f(z)| d\mu(y) d\mu(z).
\end{align*}
Using the properties of the standard kernel, we deduce that
\begin{align*}
\left|(-s\Delta) e^{s\Delta} T(f) (x_0)\right| & \leq \int\int s \left|\nabla_y  p_s(x_0,y)\right| \frac{d(y,x_0)^\epsilon}{(d(y,z)+d(x_0,z))^{d+3N+1+\epsilon}} |f(z)| d\mu(y) d\mu(z). 
\end{align*}
Then using Cauchy-Schwarz inequality and the weighted estimates on the gradient of the heat kernel (see Lemma 2.2 of \cite{CD1})
$$ \int \left|\nabla_y p_s(x_0,y)\right|^2 e^{\gamma d(x_0,y)^2/(2s)} d\mu(y) \lesssim \mu(B(x_0,r))^{-1/2}$$
we obtain (since $d(x_0,z)\geq d(Q_2,Q_1)\geq 2r$):
\begin{align*}
\left|(-s\Delta) e^{s\Delta} T(f) (x_0)\right| & \lesssim \frac{1}{\mu(B(x_0,r))} \int \frac{1}{\left(1+\frac{d(x_0,z)}{r}\right)^{d+3N+1+\epsilon}} |f(z)| d\mu(z), 
\end{align*}
which yields (\ref{assum:offfin}) by integrating over $x_0\in Q_2$.
\findem

\begin{rem} Due to our main theorem \ref{thm:principal} with $\delta>1$, we can only treat operators of order $1+\epsilon$ with $\epsilon >0$. However this result seems to be the first one, which holds in a general Riemannian manifold.
\end{rem}

\mb Usually, for a linear operator $T$ associated to a standard kernel $K$. We use the following {\it weak boundedness property}~: for all smooth functions $f,g$ such that for some $x_0\in M$ and $R>0$ and every $\alpha$
$$ \left|\nabla^\alpha f(x)\right| \lesssim \mu(B(x_0,R))^{-d-|\alpha|} \left(1+\frac{d(x,x_0)}{R}\right)^{-N}$$
for large enough integer $N$ (and similarly for $g$), we have
$$ \left|\langle T(f),g \rangle \right| \lesssim \frac{1}{\mu(B(x_0,R))}.$$
Indeed this property implies our one (\ref{assum:weaklybounded}), since $e^{-tL}(f)$ is a smooth function at the scale $t^{1/2}$. We also recover the classical T(1) theorem in Euclidean space and extend it on a large class of Riemannian manifold. It could be interesting to make a mixture of our present study with the works of F. Nazarov, S. Treil and A. Volberg \cite{NTV, NTV2} and X. Tolsa \cite{T} in order to obtain results in Riemanian manifolds with a non doubling measure.

\gb We finish this work by asking an open question: in some situations, the semigroup $(e^{-tL})_{t>0}$ does not satisfy pointwise estimates as (\ref{eq:pointwise}) but only $L^2-L^2$ off-diagonal estimates (like Gaffney estimates). Can we expect a similar $T(1)$-theorem under just off-diagonal decays for the heat kernel ? In our proof, the pointwise bound seems to be very important in a one hand to get a Sobolev inequality (Proposition \ref{prop:sobolev}) and in a other hand to bound the maximal function (\ref{eq:maxnt}) appearing in the ``Carleson measure - argument''.

\end{document}